\documentclass{amsart}
\usepackage{tikz,hyperref,enumitem,caption}
\usepackage[utf8]{inputenc}
\usepackage[T1]{fontenc}
\usepackage{microtype}
\usepackage[charter]{mathdesign}
\usetikzlibrary{decorations.pathmorphing}
\usetikzlibrary{arrows,automata,positioning}

\newlength{\proofmargin}
  \renewenvironment{thebibliography}[1]{%
    \begin{oldthebibliography}{#1}%
      \setlength{\parskip}{0ex}%
      \setlength{\itemsep}{0.1ex}%
      \setlength{\labelwidth}{1.4cm} 
  }%
  {%
    \end{oldthebibliography}%
  }

\DeclareMathOperator{\Q}{\mathbf{Q}}
\DeclareMathOperator{\F}{\mathbf{F}}
\DeclareMathOperator{\Z}{\mathbf{Z}}
\DeclareMathOperator{\C}{\mathbf{C}}

\DeclareMathOperator{\X}{\mathcal{X}}
\DeclareMathOperator{\Y}{\mathcal{Y}}
\renewcommand{\ss}{\mbox{\scriptsize ss}}

\DeclareMathOperator{\CatX}{\mathbf{Form}^{\ss}_X}
\DeclareMathOperator{\CatV}{\mathbf{Vert}^{\ss}_X}

\DeclareMathOperator{\Corr}{\mathfrak{C}}

\DeclareMathOperator{\ad}{\mbox{\scriptsize ad}}
\DeclareMathOperator{\et}{\mbox{\scriptsize \'et}}

\DeclareMathOperator{\sh}{\mbox{\scriptsize sh}}
\DeclareMathOperator{\sC}{\mbox{\scriptsize sc}}
\DeclareMathOperator{\sm}{\mbox{\scriptsize sm}}
\DeclareMathOperator{\rig}{\mbox{\scriptsize{rig}}}
\DeclareMathOperator{\Het}{\HH^1_{\mbox{\scriptsize \'et}}}
\DeclareMathOperator{\Heti}{\HH^i_{\mbox{\scriptsize \'et}}}
\DeclareMathOperator{\sing}{\mbox{\scriptsize sing}}
\newtheorem{theorem}{Theorem}

\newtheorem{conjecture}{Conjecture}

\newtheorem{lemma}{Lemma}
\newtheorem{corollary}{Corollary}
\newtheorem{thmx}{Theorem}
\newtheorem{innercustomthm}{Theorem}
\newenvironment{customthm}[1]
  {\renewcommand\theinnercustomthm{#1}\innercustomthm}
  {\endinnercustomthm}
\newcommand\HH{\mathrm{H}}

\makeatletter
\renewenvironment{proof}[1][\proofname]%
{%
\par\pushQED{\qed}\normalfont\topsep6\p@\@plus6\p@\relax%
\begin{list}{}{\rightmargin=8pt\leftmargin=\proofmargin}%
  \item[\hskip\labelsep\bfseries#1\@addpunct{.}]\ignorespaces
}{%
\popQED\end{list}\@endpefalse%
}%
\makeatother

\title{Stable models of Hecke operators}
\author{Jan Vonk}
\date{}

\begin{document}

\begin{abstract}
We investigate the geometry of \textit{correspondences} between curves, and prove that correspondences over a non-Archimedean valued field have potentially stable reduction, generalising and strengthening results of Coleman and Liu. This yields a concrete description of the operator on the cohomology of the generic fibres arising from linearisation of the correspondence, via the weight-monodromy filtration and Picard--Lefschetz theory. We explicitly determine stable models of Hecke operators on various quaternionic Shimura curves, and prove a generalisation of the geometric theory of canonical subgroups by Goren and Kassaei.
\end{abstract}
\maketitle
\tableofcontents

\section*{Introduction}

Hecke operators arise as linearisations of \textit{correspondences} between modular curves. Classical results due to Eichler--Shimura, Mestre--Oesterl\'e and Ribet exhibit deep arithmetic information by studying reductions mod $p$ of various Hecke correspondences. These methods are usually restricted to situations where the level is divisible by $p$ at most once, as they crucially rely on the semi-stability of the canonical Katz--Mazur models of the modular curves involved. Even though precise information on semi-stable models of modular curves is known in great generality by the recent work of Weinstein \cite{Wei12}, we need additional information about the degeneracy maps to carry out similar arguments. In this paper, we use non-Archimedean analytic geometry to propose a systematic framework for the investigation of operators arising from correspondences between curves, and generalise several results in the literature. 

\par We postpone linearisation as long as possible, and start by studying correspondences as \textit{geometric objects}. We prove that correspondences between curves have potentially stable reduction, generalising and strengthening the work of Coleman \cite{Col03} and Liu \cite{Liu06}, resulting in a strong analogue of the theorem of Deligne--Mumford \cite{DM69}. We introduce \textit{skeletal semi-stable models} for correspondences, and apply them to the study of the resulting spectrum via the weight-monodromy filtration on cohomology. This suggests a systematic geometric study of spectra of Hecke operators, avoiding the use of moduli interpretations which might be absent for automorphic forms that do not arise from PEL Shimura varieties. As an illustration of our methods, we generalise the work of Goren--Kassaei \cite{GK06} on canonical subgroups, following their strategy of focusing on the underlying \textit{geometry} of the considered morphisms. 

\subsection*{Stable models of correspondences. }Coleman \cite{Col03} and Liu \cite{Liu06} prove potentially semi-stable reduction theorems for finite morphisms between curves over a general class of fields. We generalise and strengthen their results to the setting of \textit{correspondences} between curves, compactifying a suitably defined moduli stack. More precisely, let $\C:Y_1 \leftarrow X \rightarrow Y_2$ be a correspondence between smooth, proper, geometrically connected curves $X,Y_1$ and $Y_2$ over a non-Archimedean field $K$, where the maps are finite. 
\begin{customthm}{A}
After a finite separable field extension of $K$, we can find a correspondence $\mathfrak{C}$:
\begin{center}
\begin{tikzpicture}[->,>=stealth']
 \node[] (X) {$\mathcal{X}$};
 \node[yshift=-.8cm,left of=X,node distance=1.6cm]  (Y1) {$\mathcal{Y}_1$};   
 \node[yshift=-.8cm,right of=X,node distance=1.6cm] (Y2) {$\mathcal{Y}_2$};

 \path (X) edge node[anchor=south,above left]{\footnotesize $\pi_1$} (Y1);
 \path (X) edge node[anchor=south,above right]{\footnotesize $\pi_2$} (Y2);
\end{tikzpicture}
\end{center}
where $\mathcal{X}, \mathcal{Y}_1$ and $\mathcal{Y}_2$ are semi-stable models over the valuation ring of $K$ for $X,Y_1$ and $Y_2$, the morphisms $\pi_1$ and $\pi_2$ are finite, such that $\mathfrak{C}$ restricts to $\C$ on the generic fibres. There is a unique such correspondence which is minimal with respect to the relation of pairwise domination of its objects.
\end{customthm}
\par \noindent In fact, we obtain a considerably stronger statement as a corollary of Theorem \ref{ssskel}. We show that we may always find a \textit{skeletal} such correspondence $\mathfrak{C}$, which satisfies the additional constraint that the singular loci of the special fibres of $\mathcal{Y}_1$ and $\mathcal{Y}_2$ pull back to the singular locus of the special fibre of $\mathcal{X}$. Skeleta of correspondences may naturally be considered as tropical objects. The association of a skeleton to a correspondence between curves may be viewed as an analogue of the tropicalisation map
\[\mathrm{Trop}:\overline{\mathcal{M}}_{g,n}^{\ad} \ \longrightarrow \  \overline{M}_{g,n}^{\mbox{\scriptsize trop}},\]
from the Deligne--Mumford--Knudsen moduli stack of $n$-pointed curves of genus $g$ to the corresponding moduli space of tropical curves. This map is studied by Abramovich, Caporaso and Payne \cite{ACP15}, and is the retraction of the moduli stack of curves onto its canonical skeleton in the sense of Thuillier \cite{Thu07}. 

\subsection*{Weight-monodromy and Shimura curves. }Skeletal semi-stable models give us precise information about the spectral properties of various linearisations of $\C$. For $l$ a prime different from the residue characteristic of $K$, we obtain a linear map 
\[\C^*: \Heti(Y_{1,\overline{K}}, \Q_l) \rightarrow \Heti(Y_{2,\overline{K}},\Q_l),\]
defined by $\pi_{2,*}\circ \pi_1^*$. This map is \textit{strict} in the sense that it respects the weight-monodromy filtrations, and the induced map on graded pieces may be described explicitly from a skeletal semi-stable model via Picard--Lefschetz theory. N\'eron component groups $\Phi_i$ are described by the $l$-adic monodromy pairing on cohomology groups, and we may consequently also compute the induced map $\C^{*}: \Phi_1 \rightarrow \Phi_2$. 

\par Applied to the setting of quaterionic Shimura curves, this allows us to revisit many classical results such as the Eichler--Shimura relation, the work of Ribet on level lowering \cite{Rib90}, and graph algorithms due to Mestre--Oesterl\'e \cite{Mes86} and Demb\'el\'e--Voight \cite{DV13}. Theorems C and D explicitly determine skeletal semi-stable models at $\mathfrak{p}$ of various Hecke operators $T_{\mathfrak{l}}$ on quaternionic Shimura curves. This is done both for $\mathfrak{l} \neq \mathfrak{p}$ and $\mathfrak{l} = \mathfrak{p}$, and the latter case exhibits extra components not present in the stable models of the individual curves involved. Letting $X$ denote a quaternionic Shimura curve with sufficiently small level away from $\mathfrak{p}$, we obtain in section \ref{4Chap} the following pictures for the stable skeleta:

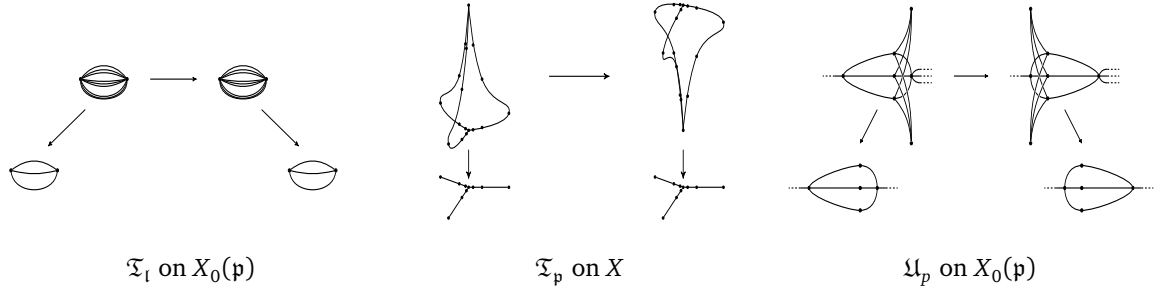
\begin{center}
\begin{figure}[!htb]
\minipage{0.32\textwidth}
\resizebox{0.9\textwidth}{0.6\textwidth}{
\begin{tikzpicture}[-,line width=0.8pt]
 \path (1,3) edge[->,>=stealth'] (3,3);
 \node[] at (2,5.5) {};
 \node[] at (2,-1.5) {};

\tikzstyle{every node}=[draw,circle,fill=black,minimum size=3pt,inner sep=0pt]

 \node[](XNL1) at (-5,0) {};
 \node[](XNL2)   at (-3,0)  {};   
 \node[](Xl1) at (-2,3)  {};
 \node[](Xl2)   at (0,3)  {};
 \node[] (Xr1) at (4,3)  {};
 \node[](Xr2)   at (6,3)  {};
 \node[](XNR1) at (7,0)  {};
 \node[](XNR2)   at (9,0)  {}; 
 
 \path (XNL1) edge[bend left=30] (XNL2) 
 	   edge[bend right=10] (XNL2) 
 	   edge[bend right=80] (XNL2);
 \path (XNR1) edge[bend left=30] (XNR2)
 	   edge[bend right=10] (XNR2) 
 	   edge[bend right=80] (XNR2);
 \path (Xl1) edge[bend left=20] (Xl2) 
 	   edge[bend left=30] (Xl2) 
 	   edge[bend left=45] (Xl2) 
 	   edge[bend right=7] (Xl2) 
 	   edge[bend right=15] (Xl2)
 	   edge[bend right=23] (Xl2)  
 	   edge[bend right=64] (Xl2) 
 	   edge[bend right=78] (Xl2) 
 	   edge[bend right=90] (Xl2);
 \path (Xr1) edge[bend left=20] (Xr2) 
 	   edge[bend left=30] (Xr2) 
 	   edge[bend left=45] (Xr2) 
 	   edge[bend right=7] (Xr2) 
 	   edge[bend right=15] (Xr2)
 	   edge[bend right=23] (Xr2)  
 	   edge[bend right=64] (Xr2) 
 	   edge[bend right=78] (Xr2) 
 	   edge[bend right=90] (Xr2);
 \path (5.8,2) edge[->,>=stealth'] (7.4,.8);
 \path (-1.8,2) edge[->,>=stealth'] (-3.4,.8);

\end{tikzpicture}}
\captionsetup{labelformat=empty}
\caption{$\mathfrak{T}_{\mathfrak{l}}$ on $X_0(\mathfrak{p})$}
\endminipage\hfill
\minipage{0.32\textwidth}
\begin{center}
\resizebox{0.8\textwidth}{0.6\textwidth}{
\begin{tikzpicture}[-,scale=0.3,node distance=1.8cm]
 \path (0,2.2)  edge[->,>=stealth'] node[left,draw=none] {} (0,0.5);
 \path (16,2.2) edge[->,>=stealth'] node[left,draw=none] {} (16,0.5);
 \path (6,6.5) edge[->,>=stealth'] (10.5,6.5);

\tikzstyle{every node}=[draw,circle,fill=black,minimum size=1.2pt,inner sep=0pt]


 \node[] (X) at (0,0)  {};
 \node[] (X1) [above of =X, node distance=1cm]  {};
 \node[] (X2) [above of =X1,node distance=2.2cm]  {};
 \node[] (Y) at (3,0)  {};   
 \node[] (P) at (1,0)  {};
 \node[] (Q) at (0.33,0)  {};
 \node[] (Y1) [above of =Y, node distance=1.3cm]  {};   
 \node[] (P1) [above of =P, node distance=1.9cm]  {};
 \node[] (P2) [above of =P, node distance=1.05cm]  {};
 \node[] (Q1) [above of =Q, node distance=2.5cm]  {};
 \node[] (Q2) [above of =Q, node distance=1.01cm]  {};
\draw plot[smooth, tension=0.9] coordinates {(X1) (Y1) (P1) (X2)};
\path (X) edge (P);
\path (P) edge (Y);

 \node[] (Y) at (-2.1,0.6)  {};   
 \node[] (P) at (-0.7,0.2)  {};
 \node[] (Q) at (-0.233,0.066)  {};
 \node[] (Y1) [above of =Y, node distance=1.3cm]  {};   
 \node[] (P1) [above of =P, node distance=1.9cm]  {};
 \node[] (P2) [above of =P, node distance=1.05cm]  {};
 \node[] (Q1) [above of =Q, node distance=2.5cm]  {};
 \node[] (Q2) [above of =Q, node distance=1.01cm]  {};
\draw plot[smooth, tension=0.9] coordinates {(X1) (Y1) (P1) (X2)};
\path (X) edge (P);
\path (P) edge (Y);

 \node[] (Y) at (-1.5,-1.8)  {};   
 \node[] (P) at (-0.5,-0.6)  {};
 \node[] (Q) at (-0.166,-.2)  {};
 \node[] (Y1) [above of =Y, node distance=1.3cm]  {};   
 \node[] (P1) [above of =P, node distance=1.9cm]  {};
 \node[] (P2) [above of =P, node distance=1.02cm]  {};
 \node[] (Q1) [above of =Q, node distance=2.5cm]  {};
 \node[] (Q2) [above of =Q, node distance=1cm]  {};
\draw plot[smooth, tension=1.3] coordinates {(X1) (Y1) (P1) (X2)};
\path (X) edge (P);
\path (P) edge (Y);


\begin{scope}[xshift=16cm]

 \node[] (X) at (0,0)  {};
 \node[] (X1) [above of =X, node distance=1cm]  {};
 \node[] (X2) [above of =X,node distance= 3.2cm]  {};
 \node[] (Y) at (3,0)  {};   
 \node[] (P) at (1,0)  {};
 \node[] (Q) at (0.33,0)  {};
 \node[] (Y1) [above of =Y, node distance=2.9cm]  {};   
 \node[] (P1) [above of =P, node distance=2.3cm]  {};
 \node[] (P2) [above of =P, node distance=3.15cm]  {};
 \node[] (Q1) [above of =Q, node distance=1.6cm]  {};
 \node[] (Q2) [above of =Q, node distance=3.185cm]  {};
 \draw plot[smooth, tension=0.9] coordinates {(X1) (P1) (Y1) (X2)};
\path (X) edge (P);
\path (P) edge (Y);

 \node[] (Y) at (-2.1,0.6)  {};   
 \node[] (P) at (-0.7,0.2)  {};
 \node[] (Q) at (-0.233,0.066)  {};
 \node[] (Y1) [above of =Y, node distance=2.9cm]  {};   
 \node[] (P1) [above of =P, node distance=2.3cm]  {};
 \node[] (P2) [above of =P, node distance=3.15cm]  {};
 \node[] (Q1) [above of =Q, node distance=1.6cm]  {};
 \node[] (Q2) [above of =Q, node distance=3.19cm]  {};
\draw plot[smooth, tension=0.9] coordinates {(X1) (P1) (Y1) (X2)};
\path (X) edge (P);
\path (P) edge (Y);

 \node[] (Y) at (-1.5,-1.8)  {};   
 \node[] (P) at (-0.5,-0.6)  {};
 \node[] (Q) at (-0.166,-.2)  {};
 \node[] (Y1) [above of =Y, node distance=2.9cm]  {};   
 \node[] (P1) [above of =P, node distance=2.3cm]  {};
 \node[] (P2) [above of =P, node distance=3.15cm]  {};
 \node[] (Q1) [above of =Q, node distance=1.6cm]  {};
 \node[] (Q2) [above of =Q, node distance=3.18cm]  {};
\draw plot[smooth, tension=1] coordinates {(X1) (P1) (Y1) (X2)};
\path (X) edge (P);
\path (P) edge (Y);

\end{scope}
\end{tikzpicture}}
\end{center}
\captionsetup{labelformat=empty}
\caption{$\mathfrak{T}_{\mathfrak{p}}$ on $X$}
\endminipage\hfill
\minipage{0.32\textwidth}%
\resizebox{\textwidth}{0.6\textwidth}{
\begin{tikzpicture}[-,line width=0.9pt]
 \path (2,3.5) edge[->,>=stealth'] node[above right,draw=none] {} (3,2);
 \path (-9,3.5) edge[->,>=stealth'] node[above left,draw=none] {} (-10,2);
 \path (-4.5,5) edge[->,>=stealth'] (-2.5,5);
 \node[] at (0,8.3) {};
 \node[] at (0,-1.3) {};

\tikzstyle{every node}=[draw,circle,fill=black,minimum size=4pt,inner sep=0pt]
 \node[] (XF) at (2,0)  {};
 \node[]   (XV) at (6,0)  {};   
 \node[]    (1) at (3,1)  {};
 \node[]    (2) at (3,0)  {};
 \node[]    (3) at (3,-1) {};
\draw plot[smooth cycle, tension=0.7] coordinates {(XF) (1) (XV) (3)};
\path (XF) edge (XV);
\path (XF) edge (1.5,0);
\path (1.5,0) edge[dashed] (0.8,0);
\path (XV) edge (6.5,0);
\path (6.5,0) edge[dashed] (7.2,0);

 \node[]  (XFF) at (0,5)  {};
 \node[]  (XVV) at (4,5)  {}; 
 \node[]  (X+) at (0,8)  {};
 \node[]  (X-) at (0,2)  {};  
 \node[]  (11) at (1,6)  {};
 \node[]  (22) at (1,5)  {};
 \node[]  (33) at (1,4)  {};
\draw plot[smooth cycle, tension=0.7] coordinates {(XFF) (11) (XVV) (33)};
\path (XVV) edge (XFF);

 \path (X+) edge[bend right=20] (11);
 \path (X+) edge[bend right=15] (22);
 \path (X+) edge[bend right=15] (33);  
 \path (X-) edge[-,line width=3pt,draw=white,bend left=15] (22);
 \path (X-) edge[-,line width=3pt,draw=white,bend left=15] (11);
 \path (X-) edge[bend left=20] (33);
 \path (X-) edge[bend left=15] (22);
 \path (X-) edge[bend left=15] (11); 
 \path (XFF) edge (-.5,5) ;
 \path (-.5,5) edge[dashed] (-1.2,5);
 
 \path (XVV) edge (4.5,5)
 			 edge[bend left] (4.5,5.3)
 			 edge[bend right] (4.5,4.7);
 \path (4.5,5) edge[dashed] (5.2,5);
 \path (4.5,5.3) edge[dashed] (5.2,5.3);
 \path (4.5,4.7) edge[dashed] (5.2,4.7);
\begin{scope}[xshift=-11cm]
 \path (0,5) edge (-.5,5);
 \path (-.5,5) edge[dashed] (-1.2,5);
 \path (4,5) edge (4.5,5)
 			 edge[bend left]  (4.5,5.3) 
 			 edge[bend right] (4.5,4.7);
 \path (4.5,5) edge[dashed]   (5.2,5);
 \path (4.5,5.3) edge[dashed] (5.2,5.3);
 \path (4.5,4.7) edge[dashed] (5.2,4.7);
\end{scope}
 
\begin{scope}[xscale=-1,xshift=7cm]
 \node[] (XF) at (2,0)  {};
 \node[] (XV) at (6,0)  {};   
 \node[] (1) at (3,1)  {};
 \node[] (2) at (3,0)  {};
 \node[] (3) at (3,-1) {};
\draw plot[smooth cycle, tension=0.7] coordinates {(XF) (1) (XV) (3)};
\path (XF) edge (XV);
\path (XF) edge (1.5,0);
\path (1.5,0) edge[dashed] (0.8,0);
\path (XV) edge (6.5,0);
\path (6.5,0) edge[dashed] (7.2,0);

 \node[]  (XFF) at (0,5)  {};
 \node[]  (XVV) at (4,5)  {}; 
 \node[]  (X+)  at (0,8)  {};
 \node[]  (X-)  at (0,2)  {};  
 \node[]  (11)  at (1,6)  {};
 \node[]  (22)  at (1,5)  {};
 \node[]  (33)  at (1,4)  {};

\draw plot[smooth cycle, tension=0.7] coordinates {(XFF) (11) (XVV) (33)};
\path (XVV) edge (XFF);

 \path (X+) edge[bend right=20] (11);
 \path (X+) edge[bend right=15] (22);
 \path (X+) edge[bend right=15] (33);  
 \path (X-) edge[-,line width=3pt,draw=white,bend left=15] (22);
 \path (X-) edge[-,line width=3pt,draw=white,bend left=15] (11);
 \path (X-) edge[bend left=20] (33);
 \path (X-) edge[bend left=15] (22);
 \path (X-) edge[bend left=15] (11);  
\end{scope} 
\end{tikzpicture}}
\captionsetup{labelformat=empty}
\caption{$\mathfrak{U}_p$ on $X_0(\mathfrak{p})$}
\endminipage
\addtocounter{figure}{-3}
\end{figure}
\end{center}

\subsection*{Overconvergent $p$-adic modular forms. }Our initial goal was a geometric study of the Atkin operator $U_p$ on spaces of \textit{overconvergent $p$-adic modular forms}. Theorem \ref{Mod} generalises a result of Coleman \cite{Col05}, and states that certain components of $\mathfrak{U}_p$ classify supersingular curves that fail to have a canonical subgroup in the sense of Katz \cite{Kat73} and Lubin \cite{Lub79}, or curves $p$-isogenous to them. These components mark the boundary of the regions $X_r$ studied in Coleman \cite{Col97} for $r=p/(p+1)$ and $r=1/(p+1)$. Consequently, we obtain an explicit description of the \textit{``boundary correspondence"} induced by $\mathfrak{U}_p$. We then generalise the work on canonical subgroups by Goren--Kassaei \cite{GK06}, and further their geometric line of investigation by proving the existence of a ``canonical subgroup section'' for general morphisms:
\begin{customthm}{F}Let $g: \mathcal{X} \rightarrow \mathcal{Y}$ be a finite flat map between semi-stable curves, and let $\mathfrak{g}: \mathcal{X}_1 \rightarrow \mathcal{Y}_1$ be the minimal skeletal semi-stable model dominating $g$. If $g$ is an isomorphism on a component $Z$ of $\mathcal{X}_{\overline{s}}$, then $g$ is an isomorphism on $\mathrm{sp}_{\mathcal{X}_1}^{-1}(Z)$. 
\end{customthm}
\par \noindent The morphism $\mathfrak{g}$ is called the \textit{skeletal stable hull} of $g$. The resulting section of $g$ is maximal in a precise sense, and the true virtue of Theorem \ref{Can} is that it identifies this maximal section to $g$ precisely, most naturally phrased in terms of the skeletal stable hull of $g$. A particular case was handled by Goren--Kassaei \cite{GK06}, where the maximal section could be identified by ad hoc methods.

\subsection*{Outline. }We recall results on the specialisation map of adic spaces and functoriality of semi-stable models of curves in section \ref{1Chap}, and prove that correspondences of curves over non-Archimedean valued fields have potentially semi-stable reduction in section \ref{2Chap}. We explain in section \ref{3Chap} how the knowledge of a skeletal semi-stable model aids us in the computation of the linear map between \'etale cohomology and N\'eron component groups arising from linearisation of a correspondence. We explicitly determine the stable models of Hecke operators on various quaternionic Shimura curves in section \ref{4Chap}, and generalise the geometric arguments of Goren--Kassaei on the existence of canonical subgroups in section \ref{5Chap}. 

\subsection*{Acknowledgements. }We are greatly indebted to Netan Dogra for countless enlightening discussions, and to Minhyong Kim for proposing a study of the geometry of correspondences. We thank Alexander Betts, Kęstutis Česnavičius, Matthew Emerton, Lorenzo Fantini and Christian Johansson for several helpful comments. The author is supported by the Engineering and Physical Sciences Research Council (EPSRC), the Andrew Mullins Scholarship, St. Catherine's College Oxford, and the Leatherseller's Company.

\section{Specialisation maps and formal fibres}\label{1Chap}

We start by recalling the notion of semi-stable vertex sets of smooth quasi-projective curves. They provide us with a combinatorial tool to understand semi-stable models and finite maps between them, via the specialisation map. Most of the material in this section is treated in more detail by Baker--Payne--Rabinoff \cite{BPR11} and Amini--Baker--Brugall\'e--Rabinoff \cite{ABBR14} using Berkovich geometry. We adopt the language of adic spaces as developed by Huber \cite{Hub93,Hub94}, as it streamlines our presentation.

\subsection{Notation. }\label{1Not}Let $K$ be an algebraically closed, complete, non-Archimedean field with topology induced by a non-trivial valuation $|\cdot |$ of rank $1$. Let $R$ be its valuation ring with maximal ideal $\mathfrak{m}$ and residue field $k$. For a curve $\mathcal{X}$ over $R$, we write $\mathcal{X}_s$ for its special fibre, and $X$ for its generic fibre. We let $f:X \rightarrow Y$ be a finite morphism between smooth proper connected curves over $K$, and $X^{\ad}, Y^{\ad}$ the adic spaces associated to $X,Y$. We will abuse notation and write $f$ for the induced map $f^{\ad}$ between these adic spaces. Given a point $x \in X^{\ad}(K)$, we write $k(x)$ for the residue field of $\mathcal{O}_{X^{\mbox{\tiny ad}},x}$ and $k(x)^+$ for the image of $\mathcal{O}_{X^{\mbox{\tiny ad}},x}^+$. 

\subsection{Hyperbolic curves and stable reduction. }\label{1Hyp}To add flexibility to our treatment of stable reduction, so as to include curves of small genus, we will allow \textit{punctures} $D_X \subset X(K)$ and $D_Y \subset Y(K)$, which are finite sets of Type-I points such that $f^{-1}(D_Y) = D_X$. Recall that a smooth, proper, punctured curve $(X,D_X)$ is said to be \textit{hyperbolic} if $\chi(X,D_X) < 0$, where $\chi(X,D_X)= 2 - 2g(X) - |D_X|$ is the Euler characteristic. A \textit{semi-stable formal model} of the punctured curve $(X,D_X)$ is an integral proper admissible formal $R$-scheme $\mathfrak{X}$ such that its generic fibre as an adic space is isomorphic to $X^{\ad}$, and moreover
\begin{itemize}
\item $\mathfrak{X}_s$ is a reduced connected curve over $k$ with at most ordinary double points for singularities, 
\item all points in $D_X$ reduce to distinct smooth points on $\mathfrak{X}_s$. 
\end{itemize}
The category $\CatX$ consists of semi-stable formal models of $(X,D_X)$, together with an isomorphism between the adic generic fibre and $X^{\ad}$. A morphism between two such models is a morphism of formal $R$-schemes that induces the identity on $X^{\ad}$. When $(X,D_X)$ is hyperbolic, there exists a unique minimal such model, which we call the \textit{stable formal model}. It is characterised by the properties of semi-stability, together with the additional constraint that the Euler characteristic of every rational component of the special fibre $\mathfrak{X}_s$, minus its singular locus and the set of punctures, is negative.

\subsection{Coleman's wide open spaces. }We recall the notion of wide open balls and annuli due to Coleman \cite{Col89}, in the language of adic spaces. A wide open ball is an adic space which is isomorphic to the complement of the set $|t| = 1$ in the open ball $\mathrm{Spa}(K\langle t \rangle, R \langle t \rangle)$. A wide open annulus is an adic space isomorphic to the complement in a wide open ball of the open set $|t| \leq p^w$ for some $w \in \mathbf{R}_{>0}$ which we call the \textit{width} of the annulus. We note that a wide open ball possesses exactly one Type-V point which is not the specialisation of any Type-II point. This Type-V point is called the \textit{apex point} of the wide open ball. Similarly, wide open annuli have exactly $2$ such points, which we also call apex points. 

\subsection{The specialisation map. }Let $\mathfrak{X}$ be an admissible formal $R$-scheme with generic fibre $X^{\ad}$. As in \cite[Theorem 2.22]{Sch12} the specialisation map is a morphism of locally ringed topological spaces
\[\mathrm{sp}_{\mathfrak{X}}: \left(X^{\ad}, \mathcal{O}_{X^{\mbox{\tiny ad}}}^+ \right) \rightarrow \left(\mathfrak{X}, \mathcal{O}_{\mathfrak{X}}\right),\]
whose fibres are called \textit{formal fibres}. When $\mathfrak{X}$ is semi-stable, it is possible to determine the nature of the formal fibres by combining the work of Bosch-L\"utkebohmert \cite[Propositions 2.2 and 2.3]{BL85} and Berkovich \cite[Proposition 2.4.4]{Ber90}. We obtain the following theorem, see also \cite[Theorem 4.6]{BPR11}.
\begin{theorem}[Bosch--L\"utkebohmert, Berkovich]\label{sp}
Let $\xi$ be a point of $\mathfrak{X}_s$. Then
\begin{itemize}
\item $\xi$ is a generic point if and only if $\mathrm{sp}^{-1}_{\mathfrak{X}}(\xi)$ consists of a single Type-II point of $X^{\ad}$,
\item $\xi$ is a smooth closed point if and only if $\mathrm{sp}^{-1}_{\mathfrak{X}}(\xi)$ is a wide open ball,
\item $\xi$ is an ordinary double point if and only if $\mathrm{sp}^{-1}_{\mathfrak{X}}(\xi)$ is a wide open annulus.
\end{itemize}
\end{theorem}

\subsection{Semi-stable vertex sets. }A \textit{semi-stable vertex set} of the smooth, proper, punctured curve $(X,D_X)$ is a finite set $V$ of Type-II points of $X^{\ad}$ such that 
\begin{itemize}
\item the space $X^{\ad} \backslash V$ is a disjoint union of wide open balls and finitely many wide open annuli,
\item the points in $D_X$ belong to distinct wide open balls in $X^{\ad} \backslash V$.
\end{itemize}
For a Type-II point $x$ in a semi-stable vertex set $V$, call its \textit{valency} the number of apex points of wide open annuli in $X^{\ad}\backslash V$ in the topological closure of $x$. The category $\CatV$ consists of semi-stable vertex sets of $(X,D_X)$, where a morphism between two such sets is given by inclusion. Theorem \ref{sp} allows us to attach to a semi-stable formal model $\mathfrak{X}$ the finite set $V_{\mathfrak{X}} := \{\mathrm{sp}^{-1}_{\mathfrak{X}}(\xi)\}_{\xi}$, where $\xi$ ranges over the generic points of the irreducible components of $\mathfrak{X}_s$. It follows immediately from Theorem \ref{sp} that $V_{\mathfrak{X}}$ is a semi-stable vertex set for $(X,D_X)$. This defines a functor between $\CatX$ and $\CatV$. It turns out that it is in fact an anti-equivalence of categories, which is proved in \cite[Theorem 4.11]{BPR11}.
\begin{theorem}The functor $\CatX  \rightarrow \CatV\ : \ \mathfrak{X} \mapsto V_{\mathfrak{X}}$ induces an anti-equivalence of categories.
\end{theorem}
Choose semi-stable formal models $\mathfrak{X}, \mathfrak{Y}$ of $(X,D_X)$ and $(Y,D_Y)$ respectively. This yields a rational map $\mathfrak{X} \dashrightarrow \mathfrak{Y}$ induced by $f$. We will investigate when this rational map extends to a morphism $\mathfrak{X} \rightarrow \mathfrak{Y}$, and whether we can arrange for such an extension to be a finite map. The equivalence we just established gives us a direct way to investigate this problem, as semi-stable vertex sets naturally live in $X^{\ad}$, where we understand the behaviour of $f$. The following theorem is proved in \cite[Theorem 5.13]{ABBR14}.
\begin{theorem}\label{abbr}
Let $\mathfrak{X},\mathfrak{Y}$ be semi-stable formal models of $(X,D_X)$ and $(Y,D_Y)$, then $f$ extends to a morphism $\mathfrak{X} \rightarrow \mathfrak{Y}$ if and only if $f^{-1}(V_{\mathfrak{Y}}) \subseteq V_{\mathfrak{X}}$. This extension is finite if and only if $f^{-1}(V_{\mathfrak{Y}}) = V_{\mathfrak{X}}$.
\end{theorem} 

\subsection{Faithfully flat descent. }\label{1Des}We will be interested in base fields $K_0$ that are not necessarily algebraically closed or complete. We now recall some relevant parts of descent theory, of which a nice exposition and more details may be found in \cite[Section 5]{ABBR14}. Let $(X,D)$ be a smooth, projective, geometrically connected curve over a field $K_0$ with a non-trivial non-Archimedean valuation of rank $1$ and valuation ring $R_0$. Let $K = \overline{K}^{\wedge}_0$ be the completion of an algebraic closure of $K_0$, which is itself algebraically closed, and let $R$ be the valuation ring of $K$. The tools outlined above can be used to analyse semi-stable models of $(X_{K},D_K)$, after which we pass back to $(X,D)$ by the theory of faithfully flat descent. More precisely, any semi-stable model of $X_K$ is isomorphic to $\mathcal{X}_R$, where $\mathcal{X}$ is a semi-stable model of $X_{K_1}$ over some finite separable extension $K_1$ of $K_0$. Moreover, the following Lemma is proved in \cite[Lemma 5.5]{ABBR14}.

\begin{lemma}\label{Desc}
Let $f: X \rightarrow Y$ be a finite morphism between smooth, projective, geometrically connected curves over $K_0$, with semi-stable $R_0$-models $\mathcal{X}$ and $\mathcal{Y}$ respectively. Suppose that $f$ extends to a finite morphism $\mathcal{X}_R \rightarrow \mathcal{Y}_R$, then $f$ also extends uniquely to a finite morphism $\mathcal{X} \rightarrow \mathcal{Y}$ defined over $R_0$. 
\begin{proof}
The rational map $f: \mathcal{X} \dashrightarrow \mathcal{Y}$ extends to a morphism if and only if the projection $\pi:\Gamma_f \rightarrow \mathcal{X}$ of its graph $\Gamma_f$ in $\mathcal{X}\times_{R_0} \mathcal{Y}$ is an isomorphism. The isomorphism $\pi_R: \Gamma_{f_R} \simeq \Gamma_f \times \mathrm{Spec}(R) \rightarrow \mathcal{X}_R$ descends to an isomorphism $\pi: \Gamma_f \rightarrow \mathrm{Spec}(R_0)$, and hence $f$ extends to a morphism $\mathcal{X}\rightarrow \mathcal{Y}$. This extension is necessarily finite, as it becomes finite after faithfully flat base change.
\end{proof}
\end{lemma}

\section{Stable models of correspondences}\label{2Chap}

In this section, we prove an analogue for correspondences of the stable reduction theorem of Deligne--Mumford \cite{DM69}. This generalises the results for finite maps proved by Coleman \cite{Col03} and Liu \cite{Liu06}. Our proof uses properties of the specialisation map from analytic geometry outlined above.

\subsection{Main definitions. }Let $K$ be a field equipped with a non-trivial non-Archimedean valuation $|\cdot |$ of rank $1$, whose valuation ring $R$ has maximal ideal $\mathfrak{m}$. A \textit{punctured correspondence} is a diagram
\begin{center}
\begin{tikzpicture}[->,>=stealth']

 \node[] (X) {\small $(X,D_X)$};
 \node[yshift=-.6cm,left of=X,node distance=2cm]  (Y1) {\small $(Y_1,D_1)$};   
 \node[yshift=-.6cm,right of=X,node distance=2cm] (Y2) {\small $(Y_2,D_2)$};
 \node[xshift=-3.2cm,yshift=-.3cm] (CC) {$\C:$};

 \path (X) edge node[above left]{\small $\pi_1$} (Y1);
 \path (X) edge node[above right]{\small $\pi_2$} (Y2);
\end{tikzpicture}
\end{center}
where
\begin{itemize}
\item $X,Y_1,Y_2$ are smooth, projective, geometrically connected $K$-curves,
\item $\pi_1, \pi_2$ are finite $K$-morphisms,
\item $D_X,D_1, D_2$ are finite sets of punctures with $\pi_1^{-1}(D_1) = D_X = \pi_2^{-1}(D_2)$.
\end{itemize}
A punctured correspondence $\C$ is said to be \textit{hyperbolic} if its objects are hyperbolic punctured curves in the sense of section \ref{1Hyp}. A \textit{semi-stable $R$-model} of the punctured correspondence $\C$ is a diagram
\begin{center}
\begin{tikzpicture}[->,>=stealth']

 \node[] (X) {$\mathcal{X}$};
 \node[yshift=-.6cm,left of=X,node distance=1.8cm]  (Y1) {$\mathcal{Y}_1$};   
 \node[yshift=-.6cm,right of=X,node distance=1.8cm] (Y2) {$\mathcal{Y}_2$};
 \node[xshift=-2.8cm,yshift=-.3cm] (CC) {$\mathfrak{C}:$};

 \path (X) edge node[anchor=south,above]{} (Y1);
 \path (X) edge node[anchor=south,above]{} (Y2);
\end{tikzpicture}
\end{center}where 
\begin{itemize}
\item $\X, \Y_1$ and $\Y_2$ are integral flat proper semi-stable $R$-curves, together with isomorphisms $\X_K \simeq X$ as well as $\Y_{1,K} \simeq Y_1$ and $\Y_{2,K} \simeq Y_2$, so that their formal completions along the special fibre are semi-stable formal $R$-models for $(X,D_X), (Y_1,D_1)$ and $(Y_2,D_2)$ in the sense of \ref{1Hyp},
\item the morphisms are \textbf{finite} and restrict to $\pi_1,\pi_2$ on the generic fibres, via the given isomorphisms. 
\end{itemize}
We say that a semi-stable $R$-model $\mathfrak{C}_1$ \textit{dominates} another semi-stable $R$-model $\mathfrak{C}_2$, if the objects of $\Corr_1$ dominate the objects of $\Corr_2$ pairwise. A semi-stable $R$-model is called \textit{stable} if it is minimal with respect to the relation of domination. The stable model of a hyperbolic correspondence is unique up to isomorphism if it exists, but in general the objects of the stable model are \textbf{not} the stable models of its objects.

\subsection{Stable models of Galois morphisms. }Before coming to a proof of potentially stable reduction of correspondences, which is the content of Theorem \ref{ss}, we prove a lemma about Galois maps $f:X \rightarrow Y$ to which the general case will be reduced. The weaker statement, not insisting that the extension of $f$ should be finite, was proved by Liu--Lorenzini \cite[Proposition 4.4]{LL99}. A proof for $K = \C_p$ is given in Coleman \cite[Section 3]{Col03}, and the proof we include is close in spirit to the paper \cite{ABBR14}. 
\begin{lemma}\label{galstab}
Let $f:(X,D_X) \rightarrow (Y,D_Y)$ be a finite Galois morphism between smooth, projective, geometrically connected, hyperbolic punctured curves over $K$. Assume $(X,D_X)$ has stable model $\mathcal{X}$ over $R$. Then there exists a unique semi-stable model $\mathcal{Y}$ of $(Y,D_Y)$ such that $f$ extends to a finite morphism $\mathcal{X} \rightarrow \mathcal{Y}$.
\begin{proof}
Extend scalars to $\overline{K}^{\wedge}$. Let $V \subset X^{\ad}$ be the stable vertex set of $(X,D_X)$, and set $W = f(V)$. There is a minimal semi-stable vertex set $W'\subset Y^{\ad}$ for $(Y,D_Y)$ containing $W$, which may be obtained by adding a finite number of Type-II points if necessary, see \cite[Lemma 3.15]{ABBR14}. We will prove that $W' = W$. Pick any element $y \in W' \backslash W$, and any element $x \in f^{-1}(y)$. We see that this means that $y$ must be of valency at least $3$, whereas $x$ must be of valency $2$. 
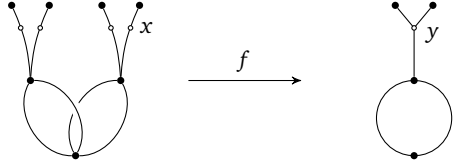
\begin{center}
\vspace{.3cm}
\begin{tikzpicture}[scale=.5]
 \node[fill,circle,scale=0.3] (O) at (0,0) {};
 \node[fill,circle,scale=0.3] (L) at (-1.2,2){};   
 \node[fill,circle,scale=0.3] (R) at (1.2,2) {};
 \node[fill,circle,scale=0.3] (LL) at (-1.7,4) {};
 \node[fill,circle,scale=0.3] (LR) at (-.7,4) {};
 \node[fill,circle,scale=0.3] (RR) at (1.7,4){};
 \node[fill,circle,scale=0.3] (RL) at ( .7,4){};

 \path (O) edge[bend left=60] (R) edge[bend right=60] (R);
 \path (O) edge[bend left=60] (L) edge[-,line width=4pt,draw=white, bend right=60] (L) edge[bend right=60] (L);
 \path (R) edge[bend left=10] (RR) edge[bend right=10] (RL);
 \path (L) edge[bend left=10] (LR) edge[bend right=10] (LL);
 
 \node[fill,circle,scale=0.3] (X) at (9,0) {};
 \node[fill,circle,scale=0.3] (Y) at (9,2) {};
 \draw (9,1) circle (1cm);
 \node[] (xlabel) at (9.5,3.2){\small $y$};
 \node[fill,circle,scale=0.3] (Z1) at (9.5,4){};
 \node[fill,circle,scale=0.3] (Z2) at (8.5,4){};
 \path (9,3.4) edge (Y) edge (Z1) edge (Z2);
 
 \path (3,2) edge[->,>=stealth'] node[above]{\small $f$} (6,2);

 \fill[white] (9,3.4) circle (.06cm);
 \draw (9,3.4) circle (.06cm);
 \fill[white] (1.46,3.4) circle (.06cm);
 \draw (1.46,3.4) circle (.06cm);
 \node[] (xtildelabel) at (1.9,3.4){\small $x$};
 \fill[white] (.94,3.4) circle (.06cm);
 \draw (.94,3.4) circle (.06cm);
 \fill[white] (-1.46,3.4) circle (.06cm);
 \draw (-1.46,3.4) circle (.06cm);
 \fill[white] (-.94,3.4) circle (.06cm);
 \draw (-.94,3.4) circle (.06cm);
\end{tikzpicture}
\vspace{.3cm}
\end{center}
Because the Galois action on fibres over Type-I points is transitive and $Y^{\ad}$ is connected, the map $f:X^{\ad} \rightarrow Y^{\ad}$ is a Galois covering of topological spaces. It follows that the valency of $x$ must be at least the valency of $y$, which is a contradiction. Therefore $W$ is semi-stable, and Theorem \ref{abbr} assures the existence of a finite morphism $f: \mathfrak{X} \rightarrow \mathfrak{Y}$ of semi-stable formal models over the ring of integral elements in $\overline{K}^{\wedge}$, where $\mathfrak{X}$ is in fact stable. We obtain algebraic models $\mathcal{X}$ and $\mathcal{Y}$ by the algebraisation theorem \cite[Corollaire 2.3.19]{Abb10}, and both $\mathcal{X}$ and $\mathcal{Y}$ descend to the integral closure of $R$ in some finite separable extension of $K$, see \cite[Lemma 5.4]{ABBR14}. Finally, the morphism $f$ descends to a finite morphism over the same field by faithfully flat descent as in Lemma \ref{Desc}.
\end{proof} 
\end{lemma}
\subsection{Remark. }We note that this theorem is false without the assumption that $f$ is Galois. However, it is true that a general $f$ extends to a morphism between the stable models of its source and target, although such an extension will in general fail to be finite. See \cite[Proposition 4.4]{LL99} for more details.

\subsection{Potentially stable reduction for correspondences. }We now come to the main result of this paper, which should be viewed as an analogue of the theorem of Deligne--Mumford \cite[Corollary 2.7]{DM69} on potentially semi-stable reduction of smooth projective curves. We note that we already made essential use of the result of Deligne--Mumford in the proof of Lemma \ref{galstab}, so we do not obtain a new proof by specialising to the degenerate case where our correspondence is a single curve with identity morphism.
\begin{thmx}\label{ss}Let $\C$ be a hyperbolic punctured correspondence over $K$. There is a finite separable extension of $K$ over which $\C$ has a stable $R$-model, which is unique up to isomorphism.
\begin{proof}Change scalars to $\overline{K}^{\wedge}$. Consider the Galois closure $g: (\widetilde{X},\widetilde{D}) \rightarrow (X,D_X)$ of both $\pi_1$ and $\pi_2$, with respect to some embedding of their function fields into $\overline{K(t)}$. This yields a diagram
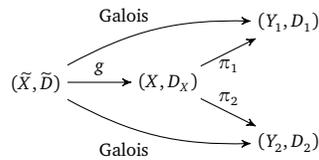
\begin{center}
\vspace{.3cm}
\begin{tikzpicture}[->,>=stealth',transform shape, scale=.8]
\node[] (Xhat) {\small $(\widetilde{X},\widetilde{D})$};
\node[right of=Xhat,node distance=2.2cm] (X) {\small $(X,D_X)$};
\node[yshift=1cm,right of=X,node distance=2cm](Y1){\small $(Y_1,D_1)$};   
\node[yshift=-1cm,right of=X,node distance=2cm](Y2){\small $(Y_2,D_2)$};
 
 \path (Xhat) edge node[anchor=south,above]{\small $g$} (X);
 \path (X) edge node[below]{{\small $\pi_1$}} (Y1);
 \path (X) edge node[above]{{\small $\pi_2$}} (Y2);
 \path (Xhat) edge[bend left=15] node[above left]{{\small Galois}} (Y1);
 \path (Xhat) edge[bend right=15] node[below left]{{\small Galois}} (Y2);
\end{tikzpicture}
\vspace{.3cm}
\end{center}
where we set $\widetilde{D} = g^{-1}(D_X)$. We see that $(\widetilde{X},\widetilde{D})$ is hyperbolic, and as such there is a unique stable vertex set $\widetilde{V} \subset \widetilde{X}^{\ad}$ by Deligne--Mumford \cite[Corollary 2.7]{DM69}. Set $V := g(\widetilde{V})$ and $W_i = \pi_i(V)$, which are semi-stable vertex sets by Lemma \ref{galstab}. We must have $V = \pi_i^{-1}(W_i)$ by transitivity of the Galois action on fibres, and hence $\pi_i$ extends to a finite morphism between the corresponding formal models. By the algebraisation theorem \cite[Corollaire 2.3.19]{Abb10}, this yields three semi-stable curves over the  ring of integral elements in $\overline{K}^{\wedge}$, which descend to a finite separable extension of $K$ by \cite[Lemma 5.4]{ABBR14}. We use faithfully flat descent as in Lemma \ref{Desc} to show that the morphisms descend to the same extension, giving us a semi-stable model $\mathfrak{C}$.

\par To show the existence of a stable model over the same field extension, we note that any semi-stable model must have the property that the semi-stable vertex sets of its objects contain the stable vertex sets. Hence, we may consider the stable vertex set $S \subset X^{\ad}$ of $(X,D_X)$ and the semi-stable vertex sets $\pi_i(S) = T_i$. Now repeat the following procedure: Enlarge $S$ to contain $\pi_i^{-1}(T_i)$, and enlarge the sets $T_i$ to contain $\pi_i(S)$. This procedure terminates, as all the newly introduced Type-II points are necessarily contained in the finite semi-stable vertex sets corresponding to the semi-stable model $\mathfrak{C}$ constructed above via the Galois closure. Therefore this procedure yields a semi-stable model $\mathfrak{C}'$ for $\C$ which is minimal with respect to the relation of domination. Hence $\mathfrak{C}'$ is the stable model of $\C$.
\end{proof}
\end{thmx}

\subsection{Remark. }We note that in general, the stable model of $\C$ does \textbf{not} consist of the stable models of its objects. As can be seen below for the case of Hecke operators, typically some extra components appear. We also note that correspondences which fail to be hyperbolic still have potentially semi-stable reduction, but in general there will not be a minimal such model. To overcome this, we can fix semi-stable models for all the objects of $\C$, and run through the above procedure to construct their \textit{stable hull} $\mathfrak{C}$, which is the minimal semi-stable model such that all the objects dominate their corresponding fixed semi-stable models. A similar relative result is proved in Theorem \ref{ssskel}.

\section{Weight-monodromy and skeleta of correspondences}\label{3Chap}

We now investigate how the knowledge of a semi-stable model for a correspondence $\C$ allows us to deduce information about the spectrum of $\C$. If $\mathfrak{C}$ is a semi-stable model of $\C$ then, even if its morphisms fail to be flat, we may use the normality of the semi-stable models involved to construct functorial trace morphisms by localisation at height $1$ primes and reduction to the flat case. Letting $f,g$ be the defining morphisms of $\mathcal{Y}_1,\mathcal{Y}_2$ to $\mathrm{Spec}(R)$, then we obtain a functorial linearisation
\[ \mathfrak{C}^*: \mathrm{R}^1f_*\Q_l \rightarrow \mathrm{R}^1g_*\Q_l,\]
which on $\overline{K}$ restricts to the map $\C^*$ on \'etale cohomology. In the case of Hecke operators, which we treat explicitly in section \ref{4Chap}, we have $Y_1 = Y_2$ and the eigenvalues are Fourier coefficients of weight two modular forms. We recall the weight-monodromy filtration on the cohomology groups of the generic fibres, and the connection with semi-stable models. In our setting, we deduce information on $\C^*$ by investigating the properties of the geometric special fibre of $\mathfrak{C}^*$, which can be explicitly understood in terms of combinatorial data coming from the dual graphs of the special fibres, and the components appearing in the reductions. The relevant notion is that of a \textit{skeletal semi-stable model} of a correspondence. 

\par \textbf{Non-Noetherian models. }As the example of the Hecke operator $T_{\mathfrak{p}}$ on quaternionic Shimura curves at $\mathfrak{p}$ in \ref{p} shows, we know that in spite of the existence of a stable model after finite separable base change as in Theorem \ref{ss}, we are often forced both to pass to an infinite field extension and allow infinitely many components when we try to construct a \textit{skeletal} such model. It might be natural to work in a larger category instead, and we hope to generalise this material in the future to more general types of adic spaces, so as to also include such perfectoid spaces as Lubin--Tate towers \cite{Wei12}.

\subsection{The weight-monodromy filtration. }\label{3WM}We now recall the definition of the monodromy filtration on the $l$-adic \'etale cohomology of a proper, smooth, separated scheme $X$ of finite type over a non-Archimedean local field $K$, and its relation with quantities that can be computed from a semi-stable model in the case of curves. By the monodromy theorem of Grothendieck \cite{SGAVII_1}, there exists a nilpotent operator 
\[N \in \mathrm{End}(\Heti(X_{\overline{K}},\Q_l))\]
such that every $\sigma$ in a sufficiently small open subgroup of the inertia group $I \subset \mathrm{Gal}(K^s/K)$ acts as $\exp(t_l(\sigma)N)$, where $t_l: I \rightarrow \Z_l(1)$ is the maximal pro-$l$ quotient map. We obtain an ascending \textit{monodromy filtration} $M_{\bullet}$ on $\Heti(X_{\overline{K}}, \Q_l)$, characterised by $NM_{i} \subseteq M_{i-2}(-1)$ and
\[ N^i: \mathrm{Gr}^M_{i} \stackrel{\sim}{\longrightarrow} \mathrm{Gr}^M_{-i}(-i) \]
From the description of the action of inertia in terms of the operator $N$, we obtain a well-defined action of the geometric Frobenius $\mathrm{Fr}_q$ determined by unipotent inertia, and in particular a notion of \textit{weights}.
\begin{conjecture}[Weight-monodromy conjecture]
The eigenvalues of $\mathrm{Fr}_q$ on $\mathrm{Gr}_j^M \Heti(X_{\overline{K}},\Q_l)$ are algebraic numbers whose conjugates all have complex absolute values equal to $q^{(i+j)/2}$.
\end{conjecture}
Let $\mathcal{X}$ be a semi-stable model for $X$ with components $Y_i$ in the special fibre, then the analysis of the vanishing and nearby cycles functor by Rapoport--Zink \cite{RZ82} yields the construction of the \textit{weight spectral sequence}, with first page
\[ E_1^{p,q} = \bigoplus_{i \geq \max(0,-p)} \HH^{q-2i}_{\et}\left(Y^{[p+2i+1]}_{\overline{s}},\Q_l(-i)\right) \hspace{.5cm} \Rightarrow \hspace{.5cm} \HH^{p+q}_{\et}\left(X_{\overline{K}},\Q_l\right),\]
where $Y^{[m]}$ is the disjoint union of the $m$-fold intersections of the $Y_i$. This induces a filtration on cohomology which is known to coincide with the weight-monodromy filtration for curves, where everything was proved by Grothendieck \cite{SGAVII_1}. Deligne proved the statement over function fields \cite{Del71a}, and some subsequent results reduce cases in mixed characteristic to this theorem, see Scholze \cite{Sch12}.

\subsection{Spectral properties of correspondences. }\label{3Spec}Let $X$ be a smooth proper curve over $K$, then the knowledge of an explicit semi-stable model $\mathcal{X}$ for $X$ helps us compute the graded pieces of the weight filtration. Indeed, the weight spectral sequence \cite{RZ82} gives us the short exact sequence
\[ 0 \rightarrow \Het(\widetilde{\mathcal{X}}_{\overline{s}}, \Q_l) \rightarrow \Het(X_{\overline{K}},\Q_l) \rightarrow \HH^1(\Gamma_{\mathcal{X}},\Q_l)(-1) \rightarrow 0, \]
where $\Gamma_{\mathcal{X}}$ is the dual graph attached to $\mathcal{X}$, and $\widetilde{\mathcal{X}}_{\overline{s}}$ denotes the normalisation of its geometric special fibre. Let $\mathcal{J}$ be the N\'eron model of the Jacobian of $X$. The quotient of $\mathcal{J}^{ }_{\overline{s}}$ by its identity connected component $\mathcal{J}^0_{\overline{s}}$ defines the N\'eron component group $\Phi$, so we get the exact sequence
\[ 0 \rightarrow  \mathcal{J}^{0}_{\overline{s}} \rightarrow \mathcal{J}^{ }_{\overline{s}} \rightarrow \Phi \rightarrow 0.\]
Recall that the results of \cite{Ray70} imply that $\mathcal{J}^{0}_{\overline{s}} \simeq \mathrm{Pic}^0(\mathcal{X}_{\overline{s}})$, see \cite[Expos\'e IX, 12.1.11]{SGAVII_1}. This means we can further decompose  $\mathcal{J}^{0}_{\overline{s}}$ as follows
\[ 0 \rightarrow \HH^1(\Gamma_{\mathcal{X}},\Z)\otimes \mathbf{G}_m \rightarrow \mathcal{J}^{0}_{\overline{s}} \rightarrow \mathrm{Pic}^0(\widetilde{\mathcal{X}}_{\overline{s}}) \rightarrow 0.\]
Note that we need not assume that $\mathcal{X}$ is regular, see \cite[Section 2]{Rib90}. The algebraic group $\HH^1(\Gamma_{\mathcal{X}},\Z)\otimes \mathbf{G}_m$ is called the \textit{toric part} of $\mathcal{J}$. The N\'eron component group $\Phi$ is isomorphic to the cokernel of the map from the toric part to its dual, given by the monodromy pairing \cite[Th\'eor\`eme 11.5]{SGAVII_1}. 

\par We see that an explicit knowledge of a semi-stable model $\mathfrak{C}$ for $\C$ is exactly what is needed to describe the induced map on the first step of the weight-monodromy filtration. The reduction of $\mathfrak{C}$ decomposes as a finite set of correspondences between irreducible curves over the residue field of $R$, providing a geometric description of the induced morphism $\mathrm{Pic}^0(\widetilde{\mathcal{Y}}_{1,\overline{s}}) \rightarrow \mathrm{Pic}^0(\widetilde{\mathcal{Y}}_{2,\overline{s}})$. For the case of the Hecke operator $T_{\mathfrak{p}}$ on modular curves $X_0^B$ treated in \ref{p}, this yields the celebrated Eichler--Shimura relation $T_{\mathfrak{p}} = \mathrm{Frob} + \mathrm{Ver}$. To describe the induced morphism between both the graded pieces of weight $2$ and the N\'eron component groups $\Phi_1 \rightarrow \Phi_2$, we use Raynaud's theorem and the monodromy pairing to reduce the problem to finding a geometric interpretation of the maps between the toric parts of the Jacobians. As recalled above, we may describe the toric part in terms of the dual graph of a semi-stable model, which has a natural interpretation as a \textit{skeleton} in the setting of adic spaces. We now recall the relevant theory, and use it to motivate the definition of a \textit{skeletal semi-stable model} in \ref{3Skel}. 

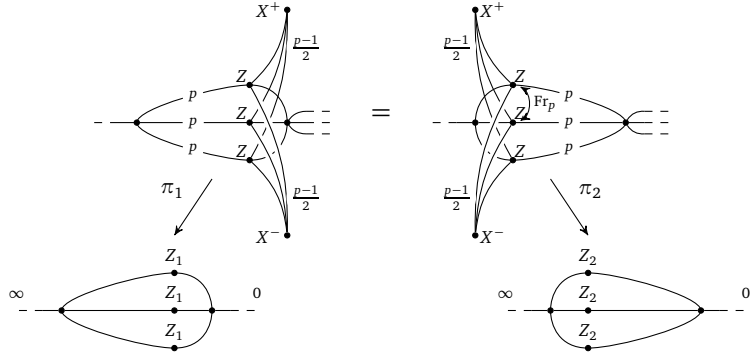
\subsection{The skeleton of a semi-stable vertex set. }To a semi-stable vertex set $V$ for $(X,D_X)$, we can associate a combinatorial structure which we call the \textit{skeleton} $\Sigma_{V}$ of $X$ with respect to $V$. First, we define $\Gamma_V$ as the set of points of $X^{\ad}$ that are not contained in a wide open ball which is disjoint from $V \cup D_X$, endowed with the subspace topology. Recall that the maximal Hausdorff quotient $\overline{\Gamma}_V$, often referred to as the \textit{Berkovich skeleton} of $X$ with respect to $V$, can naturally be viewed as a metric realisation of the dual graph. Its vertex set is $V \cup D_X$, and edges come in two flavours: There is an edge $e_Q$ between two vertices in $V$ for every intersection point $Q$ of the corresponding components in the special fibre of $\mathfrak{X}_V$, and we set the length $l(e_Q)$ of $e_Q$ to be the width of the wide open annulus $\mathrm{sp}^{-1}_V(Q)$. Every vertex in $D_X$ is adjacent to exactly one edge of length $\infty$. Define the skeleton of $X^{\ad}$ with respect to $V$, or $\mathfrak{X}_V$, to be the pair $\Sigma_V = (\Gamma_V,\mathcal{L}_V)$, where $\mathcal{L}_V = \{l(e) : e\ \mbox{edge of } \overline{\Gamma}_V\}$ consists of the set of edge lengths appearing in the maximal Hausdorff quotient $\overline{\Gamma}_V$, viewed as a metric graph. We define the \textit{genus} of a Type-II point $x$ in $\Gamma_{V}$ to be the genus of the residue field $k(x)$ as defined in \ref{1Not}. 

\subsection{Example. }Consider the projective closure of $y^2 = x^3 + x^2 + p^3$, over $\mathcal{O}_{\C_p}$, and let $\mathcal{X}$ be the blow-up at a smooth point of the special fibre. Then $\mathcal{X}$ is a normal semi-stable model whose special fibre consists of a nodal curve $C_1$ and a projective line $C_2$, crossing transversally. Now set $X$ to be the generic fibre of $\mathcal{X}$, and $D = \{ (0,1,0) \}$ the point at infinity. The formal completion $\mathfrak{X}$ of $\mathcal{X}$ is a semi-stable formal model for $(X,D)$, and its corresponding skeleton $\Sigma_{\mathfrak{X}}$ can be visualised as

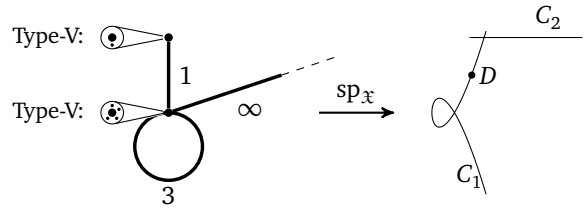
\begin{center}
\begin{tikzpicture}[scale=0.5]
  
  \draw (0,2) node[fill,circle,scale=0.35] (Y) {};
  \draw (-1.42, 2.29) edge node[]{} (-.15,2.04);
  \draw (-1.42, 1.71) edge node[]{} (-.15,1.96);
  \draw (-3.2,2) node () {\small $\mbox{Type-V:}$};
  \draw (-1.5,2) circle (0.3cm);
  \draw (-1.5,2) node[fill,circle,scale=0.35] () {};
  \draw (-1.5,1.8) node[fill,circle,scale=0.14] () {};
  
  \draw (0,0) node[fill,circle,scale=0.35] (X) {};
  \draw (0,-2.2) node (bluh) {$3$};
  \draw[line width=1.3pt] (X) edge node[anchor=east,right]{$1$} (Y);
  \draw[line width=1.3pt] (0,-0.9) circle (0.9cm);
  \draw[line width=1.3pt] (X) edge node[below right]{$\infty$} (3,1);
  \draw[dashed] (3,1) edge node[below right]{} (4.5,1.5);
  \draw (-1.42,-0.29) edge[line width=1.4pt,draw=white] node[]{} (-.15,0);
  \draw (-1.42,-0.29) edge node[]{} (-.15,-.02);
  \draw (-1.42, 0.29) edge node[]{} (-.15,0.04);
  \draw (-3.2,0) node () {\small $\mbox{Type-V:}$};
  \draw (-1.5,0) circle (0.3cm);
  \draw (-1.5,0) node[fill,circle,scale=0.35] () {};
  \draw (-1.5,0.19) node[fill,circle,scale=0.14] () {};
  \draw (-1.32,0.09) node[fill,circle,scale=0.14] () {};
  \draw (-1.62,-.15) node[fill,circle,scale=0.14] () {};
  \draw (-1.38,-.15) node[fill,circle,scale=0.14] () {}; 
  
  \draw (4,0) edge[->,thick,>=stealth'] node[above] {$\mathrm{sp}_{\mathfrak{X}}$} (6,0);
  
  \draw [domain=0:1.55, samples=40] plot ({2+0.6*\x^2+5}, {\x*(\x^2-1)} );
  \draw [domain=0:1.55, samples=40] plot ({2+0.6*\x^2+5}, {-\x*(\x^2-1)} );
  \draw (8,2) edge[-] node[above right] {$C_2$} (11,2);
  \draw (8,-1.7) node () {$C_1$};
  \draw (8.5,1) node () {$D$};
  \draw (8.06,1) node[fill,circle,scale=0.3] () {};
\end{tikzpicture}
\end{center}

\subsection{Harmonic morphisms. }\label{3Harm}Given a finite morphism $f:\mathfrak{X} \rightarrow \mathfrak{Y}$ of semi-stable formal models for $(X,D_X)$ and $(Y,D_Y)$, we know that $f^{-1}(V_{\mathfrak{Y}}) = V_{\mathfrak{X}}$. In general, it is not true that $f^{-1}(\Gamma_{\mathfrak{Y}}) = \Gamma_{\mathfrak{X}}$, see \cite[Remark 5.23]{ABBR14}. When this holds, we say that $f$ is a \textit{skeletal} finite morphism. For maps $f$, we can always modify $f$ so as to make it skeletal, see \cite[Theorem 2.4]{CKK15} and \cite[Corollary 4.18]{ABBR14}. Let $x$ be an apex point of a wide open annulus in $X^{\ad}\backslash V_{\mathfrak{X}}$, and set $l_x$ to be the width of the annulus. It can be shown that the ratio $d_x := l_{f(x)}/l_x$ equals the ramification index of the map $f:\mathfrak{X}_s \rightarrow \mathfrak{Y}_s$ at the point $\mathrm{sp}_{\mathfrak{X}}(x)$, see \cite[Lemma 2.1]{Col05}. We use this to define $d_x$ for all apex points attached to $\mathfrak{X}$, by simply setting $d_x$ to be the corresponding ramification index for the map $f$. The integer $d_x$ is often referred to as the \textit{expansion factor} of the edge of $\Gamma_{\mathfrak{X}}$ adjacent to $x$. It can be shown that the morphism $f: \Gamma_{\mathfrak{X}} \rightarrow \Gamma_{\mathfrak{Y}}$ is \textit{harmonic} in the sense that the sum
\[ \sum_{x \mapsto y} d_x,\] 
is independent of the Type-V point $y \in \overline{V}_{\mathfrak{Y}}\cap \Gamma_Y$, and equals the degree of $f$. A skeletal finite morphism $f:\mathfrak{X} \rightarrow \mathfrak{Y}$ gives rise to a harmonic morphism on skeleta in this sense, see \cite[Section 4.27]{ABBR14}. 
\subsection{Remark. }The related notion of harmonic morphisms of metrised complexes of curves is discussed in \cite[Section 2.16]{ABBR14}. The list of properties included in their definition of this abstract category is considerably longer, as they aim to reconstruct the curve from the tropical object. We shall not be concerned with such lifting questions, as our metrised complexes always arise from algebraic curves. 

\subsection{Skeletal stable models of correspondences. }\label{3Skel}A semi-stable model $\mathfrak{C}$ of a punctured correspondence $\C$ is \textit{skeletal} if $\pi_1$ and $\pi_2$ are skeletal in the sense of \ref{3Harm}. Not every semi-stable model $\mathfrak{C}$ of $\C$ is skeletal,and the question arises whether it is always possible to modify $\mathfrak{C}$ to make it skeletal, and whether there is always a minimal such model when $\C$ is hyperbolic. This was proved for finite maps by Cornelissen--Kato--Kool \cite[Theorem 2.4]{CKK15}, and Amini--Brugall\'e--Payne--Rabinoff \cite[Corollary 4.18]{ABBR14}.

\par In the proof of Theorem \ref{ss}, we found semi-stable vertex sets $V, W_1,W_2$ for the objects of $\C$, such that $\pi_1^{-1}(W_1) = \pi_2^{-1}(W_2) = V$. From \cite[Remark 4.19]{ABBR14}, we know that $\pi_i(\overline{\Gamma}_{V})$ is the union of $\overline{\Gamma}_{W_i}$ and a finite set of edges. Add the endpoints of those edges to $W_i$, and run through the following procedure: Enlarge $V$ to contain $\pi_i^{-1}(W_i)$, and enlarge the sets $W_i$ to contain $\pi_i(V)$. If this procedure terminates for a correspondence $\C$ over $K$, we obtain a skeletal semi-stable model $\mathfrak{C}$, which gives rise to a diagram
\begin{center}
\begin{tikzpicture}[->,>=stealth']

 \node[] (X) {$\Sigma_{\mathfrak{X}}$};
 \node[yshift=-.8cm,left of=X,node distance=1.6cm]  (Y1) {$\Sigma_{\mathfrak{Y}_1}$};
 \node[yshift=-.8cm,right of=X,node distance=1.6cm] (Y2) {$\Sigma_{\mathfrak{Y}_2}$};
 \node[xshift=-2.8cm,yshift=-.3cm] (CC) {$\Sigma_{\mathfrak{C}}:$};

 \path (X) edge node[anchor=south,above]{\tiny $\sigma_1$} (Y1);
 \path (X) edge node[anchor=south,above]{\tiny $\sigma_2$} (Y2);
\end{tikzpicture}
\end{center}
where $\sigma_1,\sigma_2$ are finite harmonic morphisms. We call $\Sigma_{\mathfrak{C}}$ the \textit{skeleton} of the skeletal semi-stable model $\mathfrak{C}$. 

\par The question now becomes whether this procedure always terminates. The answer in general is \textit{no}, as the example in section \ref{p} of $T_\mathfrak{p}$ on quaternionic Shimura varieties $X^B$ shows. However, in some cases the answer is \textit{yes}. Assume that either (i) $\pi_1,\pi_2$ are Galois, or (ii) either $\pi_1$ or $\pi_2$ is the identity. Every time a new vertex $v \in V$ is introduced, it appears on the stable skeleton of $(X,D_X)$, and the set of edge lengths $\mathcal{L}_V$ of the skeleton $\Sigma_V = (\Gamma_V,\mathcal{L}_V)$ remains constant by Lemma \ref{galstab}. This shows that the procedure terminates. Starting with semi-stable models $\mathcal{X},\mathcal{Y}_1, \mathcal{Y}_2$ for the objects of $\C$, this argument is insensitive to enlarging the stable vertex sets of $X,Y_1,Y_2$ to contain the vertex sets of the chosen models at the start of the procedure. We summarise this discussion in the following statement. 

\begin{thmx}\label{ssskel}
Let $\C$ be a hyperbolic punctured correspondence over $K$, and assume that either (i) both morphisms are Galois, or (ii) one of the morphisms is the identity. Given semi-stable models $\mathcal{X}, \mathcal{Y}_1, \mathcal{Y}_2$ for the objects of $\C$, there is a finite separable extension of $K$ over which $\C$ has a \textbf{skeletal} semi-stable model whose objects dominate $\mathcal{X}, \mathcal{Y}_1, \mathcal{Y}_2$ pairwise, and which is minimal with respect to the relation of domination.
\end{thmx}

This semi-stable model is called the \textit{skeletal stable hull} of $\C$ with respect to the triple $(\mathcal{X}, \mathcal{Y}_1, \mathcal{Y}_2)$. This strengthens the relative theorem for finite maps by Liu \cite{Liu06}. When we take $(\mathcal{X}, \mathcal{Y}_1, \mathcal{Y}_2)$ to be the stable models of the objects in $\C$, we call this model the \textit{skeletal stable model}. As in the proof of Theorem \ref{ss}, we may construct skeletal semi-stable models for general hyperbolic correspondences by reducing to the Galois case. However, we make no claims about the minimality of these models.

\begin{corollary}Let $\C$ be a hyperbolic punctured correspondence over $K$. There is a finite separable extension of $K$ over which $\C$ has a \textbf{skeletal} semi-stable $R$-model. 
\begin{proof}Resume the notation of the proof of Theorem \ref{ss}, and let $(\widetilde{X},\widetilde{D}) \stackrel{g}{\longrightarrow} (X,D)$ be the Galois closure of the morphisms of $\C$. Let $\widetilde{\mathfrak{C}}$ be the skeletal stable model of the resulting Galois correspondence, which is guaranteed by Theorem \ref{ssskel}, and let $\Gamma$ be the image of the skeleton of $(\widetilde{X},\widetilde{D})$. It is now easy to see that this defines a skeleton of $(X,D)$ which, together with the skeleta for $(Y_1,D_1)$  and $(Y_2,D_2)$ implicit in $\widetilde{C}$, define a skeletal semi-stable model of $\C$. 
\end{proof}
\end{corollary}

\section{Hecke operators on Shimura curves}\label{4Chap}


\par We now compute explicit examples of stable models of Hecke correspondences, and revisit classical results such as the Eichler--Shimura relation, the Mestre--Oesterl\'e \textit{m\'ethode des graphes} \cite{Mes86}, the algorithms due to Demb\'el\'e--Voight on Fourier coefficients of Hilbert modular forms \cite{DV13}, the analysis of the action on N\'eron component groups by Ribet \cite{Rib90}, and the moduli interpretation of the components of the semi-stable model of $X_0(Np^2)$ constructed by Edixhoven \cite{Edi90} due to Coleman \cite{Col05}. 

\par The explicit determination of stable models of correspondences $\C$ between curves is a difficult problem, and already the degenerate case of a single curve has been the subject of much research in arithmetic geometry. For correspondences, we not only need to describe semi-stable models for the objects, but we also need a good understanding of the geometry of the morphisms. This task is simplified when the objects and morphisms of $\C$ have a clear moduli interpretation. For all the examples discussed below, descriptions of integral models may be found by identifying a closely related PEL Shimura variety, for which one may describe integral models by computing the deformation theory of the associated $p$-divisible groups with appropriate Drinfeld level structure using Serre--Tate theory. 

\par We start by reviewing parts of the literature on semi-stable models of Shimura curves in section \ref{MSC}. We then determine the stable models at $\mathfrak{p}$ for Hecke operators $T_{\mathfrak{l}}$ for $\mathfrak{l} \neq \mathfrak{p}$ in section \ref{ell}, and $\mathfrak{l} = \mathfrak{p}$ in section \ref{p}. As an application, we generalise a theorem of Coleman in section \ref{Moduli}. 

\subsection{Models of quaternionic Shimura curves. }\label{MSC}We recall some of the literature on integral models of quaternionic Shimura curves, due to Katz--Mazur \cite{KM85}, Carayol \cite{Car86}, Edixhoven \cite{Edi90}, Buzzard \cite{Buz97} and Jarvis \cite{Jar04}. Let us first fix some notation. Let $F$ be a totally real number field of degree $d$ over $\Q$, and $B_{/F}$ a quaternion algebra which is split at exactly one infinite prime. Let $R$ be a maximal order of $B$, and set $G := \mathrm{Res}_{F/\Q}(B^{\times})$ to be the Weil restriction of the unit group of $B$. Let $\mathbf{A}^{\infty}$ be the ring of finite adeles over $\Q$. For $K$ a compact open subgroup of $G(\mathbf{A}^{\infty})$, the quaternionic Shimura curve attached to $B$ and $K$ is an algebraic curve denoted $X^B_{K}$, whose $\C$-points are
\[ X^B_K(\C) =  G(\Q ) \backslash G(\mathbf{A}^{\infty}) \times \mathfrak{H}^{\pm} / K.\]
By the theory of canonical models due to Shimura \cite{Shi70}, the curve $X^B_K$ is defined over $F$. Henceforth, we will fix a prime $\mathfrak{p}$ above a rational prime $p$ at which $B$ is split, and a sufficiently small level structure away from $\mathfrak{p}$. We simply denote $X^B$ for the corresponding Shimura curve, and $X^B_0(\mathfrak{p})$ for the Shimura curve with additional Iwahori level structure at $\mathfrak{p}$. For convenience, we will consider integral models over the strict Henselisation $\mathcal{O}_{F,\mathfrak{p}}^{\mathrm{\sh}}$, and set $W$ to be the ring of Witt vectors of $\overline{\F}_p$.

\par When $d=1$, the curves $X^B_0(p)$ are PEL Shimura varieties defined over $\Q$. They are fine moduli schemes for elliptic curves when $B = \mathrm{M}_2(\Q)$ and for so called \textit{false elliptic curves} when $B$ is non-split. False elliptic curves are abelian schemes of relative dimension $2$ with quaternionic multiplication by $R$, together with a level structure. Once level structures are understood in the sense of Drinfeld \cite{Dri74}, we obtain integral models $\mathcal{X}^B_0(p)$ over $W$. Buzzard \cite{Buz97} notes that the $p$-divisible group of a false elliptic curve $A$ over an algebraically closed field of characteristic $p$ is of the form $E[p^{\infty}] \times E[p^{\infty}]$, where $E$ is an elliptic curve. This allows him to apply the theorem of Serre--Tate to prove the following result, which is \cite[Theorem 4.7]{Buz97}. The case where $B$ is split is due to Deligne--Rapoport \cite{DR73}.
\begin{theorem}[Deligne--Rapoport, Buzzard]
The scheme $\mathcal{X}^B_0(p)$ over $W$ is regular, with special fibre consisting of two irreducible components isomorphic to $\mathcal{X}^B_{s}$, crossing transversally at the supersingular points. The degeneration map $\pi_1:\mathcal{X}^B_0(p) \rightarrow \mathcal{X}^B$ is finite flat.
\end{theorem}
We recall that Edixhoven determined a semi-stable model for $X_0(Np^2)$, defined over $W[\varpi]$ where $\varpi^{(p^2-1)/2} = p$. The special fibre contains four ordinary components. Two of these, which we call the \textit{inner ordinary components}, are isomorphic to $\mathcal{X}_0(N)_s$. The other two, $X^+$ and $X^-$, are called the \textit{outer ordinary components}, and are isomorphic to the Igusa curves $\overline{X}_0(N,\mathrm{Ig}(p)/\pm)$. There is one \textit{supersingular component} $Z^{\sigma}$ for every supersingular point $\sigma$ on $\mathcal{X}_0(Np^2)_s$, for which equations may be found in Edixhoven \cite[Theorem 2.1.1]{Edi90}. Every supersingular component intersects every ordinary component, and there are no other intersection points. The model is regular, except where the supersingular components meet the two inner ordinary components. The completed local ring at these points is isomorphic to
\[ W[\varpi] \llbracket x,y \rrbracket / \left(xy - \varpi^{\frac{p-1}{2}}\right). \]

\par When $d>1$, the situation is more complex, as the Shimura curves $X^B_0(\mathfrak{p})$ are no longer PEL. Carayol \cite{Car86} constructs integral models for these curves by comparing them to unitary Shimura curves which arise from picking an auxiliary CM field $E$ containing $F$ in which $\mathfrak{p}$ splits. These unitary Shimura curves are PEL, and integral models may be described by imposing level structure in the sense of Drinfeld \cite{Dri74} and using Serre--Tate theory. Jarvis \cite{Jar99, Jar04} makes this analysis explicit in the case of Iwahori level structure at $\mathfrak{p}$. The following theorem may be found in \cite[Section 10]{Jar99}.
\begin{theorem}[Carayol, Jarvis]
There is a regular model $\mathcal{X}^B_0(\mathfrak{p})$ over $\mathcal{O}^{\mathrm{\sh}}_{F,\mathfrak{p}}$, with a finite flat degeneration map $\pi_1: \mathcal{X}^B_0(\mathfrak{p}) \rightarrow \mathcal{X}^B$ to the model constructed by Carayol \cite{Car86}. The special fibre of $\mathcal{X}^B_0(\mathfrak{p})$ consists of two components isomorphic to $\mathcal{X}^B_{s}$, crossing transversally at the supersingular points. 
\end{theorem}

\subsection{Hecke operators away from $\mathfrak{p}$. }\label{ell}Let $\mathfrak{l}$ be coprime to both $\mathfrak{p}$ and the implicit tame level. The analysis of the special fibre of $\mathcal{X}^B_0(\mathfrak{p})$ by Deligne--Rapoport, Buzzard and Jarvis immediately yields a skeletal semi-stable model $\mathfrak{T}_{\mathfrak{l}}$ for the Hecke operators $T_{\mathfrak{l}}$ on $X^B_0(\mathfrak{p})$.
\begin{thmx}
The correspondence $T_{\mathfrak{l}}$ on $X^B_0(\mathfrak{p})$ has a skeletal semi-stable model over $\mathcal{O}_{F,\mathfrak{p}}^{\mathrm{\sh}}$, with skeleton as depicted in Figure \ref{picTl}.
\begin{proof}
We note that the regular models for $X^B_0(\mathfrak{p})$ described above are semi-stable over $\mathcal{O}_{F,\mathfrak{p}}$, regardless of the tame level structure imposed. Let $X^B_0(\mathfrak{p},\mathfrak{l})$ be the Shimura curve obtained by imposing additional Iwahori level structure at $\mathfrak{l}$, then the map $\pi_1: X^B_0(\mathfrak{p},\mathfrak{l})\longrightarrow X^B_0(\mathfrak{p})$ is finite, and as $\pi_2$ is the Atkin-Lehner involution at $\mathfrak{p}$ composed with $\pi_1$, we obtain a semi-stable model $\mathfrak{T}_{\mathfrak{l}}$ which is stable whenever $X^B_0(\mathfrak{p})$ is. The description of the special fibre now follows from the above theorems. 
\end{proof}
\end{thmx}

\begin{figure}[htbp!]
\begin{center}
\vspace{.3cm}
\begin{tikzpicture}[scale=0.5]
\tikzstyle{every node}=[draw,circle,fill=black,minimum size=2pt,inner sep=0pt]

 \node[label={180:\footnotesize $\mathcal{X}^B_s$}](XNL1) at (-5,0) {};
 \node[label={0:\footnotesize $\mathcal{X}^B_s$}](XNL2)   at (-3,0)  {};   
 \node[label={100:\footnotesize $\mathcal{X}^B_0(\mathfrak{l})_s$}](Xl1) at (-2,3)  {};
 \node[label={80:\footnotesize $\mathcal{X}^B_0(\mathfrak{l})_s$}](Xl2)   at (0,3)  {};
 \node[label={100:\footnotesize $\mathcal{X}^B_0(\mathfrak{l})_s$}] (Xr1) at (4,3)  {};
 \node[label={80:\footnotesize $\mathcal{X}^B_0(\mathfrak{l})_s$}](Xr2)   at (6,3)  {};
 \node[label={180:\footnotesize $\mathcal{X}^B_s$}](XNR1) at (7,0)  {};
 \node[label={0:\footnotesize $\mathcal{X}^B_s$}](XNR2)   at (9,0)  {}; 
 
 \path (XNL1) edge[bend left=30] (XNL2) 
 	   edge[bend right=10] (XNL2) 
 	   edge[bend right=80] (XNL2);
 \path (XNR1) edge[bend left=30] (XNR2)
 	   edge[bend right=10] (XNR2) 
 	   edge[bend right=80] (XNR2);
 \path (Xl1) edge[bend left=20] (Xl2) 
 	   edge[bend left=30] (Xl2) 
 	   edge[bend left=45] (Xl2) 
 	   edge[bend right=7] (Xl2) 
 	   edge[bend right=15] (Xl2)
 	   edge[bend right=23] (Xl2)  
 	   edge[bend right=64] (Xl2) 
 	   edge[bend right=78] (Xl2) 
 	   edge[bend right=90] (Xl2);
 \path (Xr1) edge[bend left=20] (Xr2) 
 	   edge[bend left=30] (Xr2) 
 	   edge[bend left=45] (Xr2) 
 	   edge[bend right=7] (Xr2) 
 	   edge[bend right=15] (Xr2)
 	   edge[bend right=23] (Xr2)  
 	   edge[bend right=64] (Xr2) 
 	   edge[bend right=78] (Xr2) 
 	   edge[bend right=90] (Xr2);
 \path (5.8,2) edge[->,>=stealth'] (7.4,.8);
 \path (-1.8,2) edge[->,>=stealth'] (-3.4,.8);

 \path (1,3) edge[->,>=stealth'] node[fill=none,below=-.3cm,draw=none] {\mbox{\scriptsize Edge mixing}}(3,3);

\end{tikzpicture}
\vspace{.3cm}
\end{center}
\caption{The skeleton of $\mathfrak{T}_{\mathfrak{l}}$} \label{picTl}
\end{figure}

\par From the above picture, we may recover that $T_{\mathfrak{l}}$ acts on the N\'eron component group as multiplication by $\mathrm{Nm}_{F/\Q}(\mathfrak{l})$, generalising Ribet \cite[Theorem 3.12]{Rib90}. We also obtain a geometric description of a special case of the Jacquet--Langlands correspondence, and the graph algorithms due to Mestre--Oesterl\'e \cite{Mes86} and Demb\'el\'e--Voight \cite{DV13}. Let $\overline{B}$ be the quaternion algebra over $F$ obtained from $B$ by making the split infinite place and $\mathfrak{p}$ ramified, and set $\overline{G} = \mathrm{Res}_{F/ \Q}(\overline{B}^{\times})$. The Jacquet--Langlands correspondence now provides us with a functor from automorphic representations for $\overline{G}$ to automorphic representations for $G$ which are square integrable modulo center at $\mathfrak{p}$ and the split infinite place. We may identify the action on the toric part $\HH^1(\Gamma_{\mathcal{X}^B_0(\mathfrak{p})},\mathbf{Z})\otimes \mathbf{G}_m$ of $\mathcal{J}_s^0$ with the action on the set of supersingular points
\[ \overline{G}(\Q) \backslash \overline{G}(\mathbf{A}^{\infty}) / K^{\mathfrak{p}} \times \mathcal{O}^{\times}_{\overline{B}_{\mathfrak{p}}},\]
which is a finite set amenable to computation. In fact, elements in this set may be represented as right ideal classes in an Eichler order of $\overline{B}$ depending on the level, as is explained in \cite{DV13}. By means of an example, set $\alpha= \zeta_7+\zeta_7^{-1}$ and $F=\Q(\zeta_7)^+=\Q(\alpha)$ to be the totally real cubic number field unramified outside $7$, and let $f$ be the unique normalised Hilbert modular form of prime level $3\mathcal{O}_F$ and parallel weight $(2,2,2)$. Set $B$ to be the quaternion algebra ramified at two of the three infinite places. By computing the relevant Brandt matrices for $\overline{G}$, we find the following Fourier coefficients for $f$:
\begin{center}
\begin{table}[!htbp]
\centering
\begin{tabular}{c|cccccc}
$\mathfrak{l}$ & $(2)$ & $(5)$ & $(11)$ & $(13,3+\alpha)$ & $(17)$ & $(19)$ \\
\hline
$a_{\mathfrak{l}}$ & $-4$ & $22$ & $58$ & $1$ & $0$ & $-56$\\
\end{tabular}
\end{table}
\end{center}

\subsection{Hecke operators at $\mathfrak{p}$. }\label{p}The morphisms $\pi_1, \pi_2: \mathcal{X}^B_0(\mathfrak{p}) \rightrightarrows \mathcal{X}^B$ defining $T_{\mathfrak{p}}$ on $X^B$ are finite flat, and hence define a semi-stable model for $T_{\mathfrak{p}}$. If we attempt to make this model skeletal using the procedure outlined in \ref{3Skel}, we keep introducing components which accumulate at the ordinary components. Figure \ref{picTp} depicts the first few iterations of this process. Though sufficient for applications, this non-Noetherian model could be seen as an indication that we should work in a bigger category, see section \ref{4Inf}.
\begin{figure}[htbp!]
\begin{center}
\begin{tikzpicture}[-,scale=.3,node distance=1.8cm]
 \path (0,2.2) edge[->,>=stealth'] node[left,draw=none] {\footnotesize $\pi_1$} (0,0.5);
 \path (16,2.2) edge[->,>=stealth'] node[left,draw=none] {\footnotesize $\pi_2$} (16,0.5);
 \node[label={=}] at (8,7) {};

\tikzstyle{every node}=[draw,circle,fill=black,minimum size=1.2pt,inner sep=0pt]


 \node[label={-60:\footnotesize $\mathcal{X}^B_s$}] (X) at (0,0)  {};
 \node[label={-30:\footnotesize $\mathcal{X}^B_s$}] (X1) [above of =X, node distance=1cm]  {};
 \node[label={90:\footnotesize $\mathcal{X}^B_s$}] (X2) [above of =X1,node distance=2.2cm]  {};
 \node[] (Y) at (3,0)  {};   
 \node[] (P) at (1,0)  {};
 \node[] (Q) at (0.33,0)  {};
 \node[] (Y1) [above of =Y, node distance=1.3cm]  {};   
 \node[] (P1) [above of =P, node distance=1.9cm]  {};
 \node[] (P2) [above of =P, node distance=1.05cm]  {};
 \node[] (Q1) [above of =Q, node distance=2.5cm]  {};
 \node[] (Q2) [above of =Q, node distance=1.01cm]  {};
\draw plot[smooth, tension=0.9] coordinates {(X1) (Y1) (P1) (X2)};
\path (X) edge (P);
\path (P) edge (Y);

 \node[] (Y) at (-2.1,0.6)  {};   
 \node[] (P) at (-0.7,0.2)  {};
 \node[] (Q) at (-0.233,0.066)  {};
 \node[] (Y1) [above of =Y, node distance=1.3cm]  {};   
 \node[] (P1) [above of =P, node distance=1.9cm]  {};
 \node[] (P2) [above of =P, node distance=1.05cm]  {};
 \node[] (Q1) [above of =Q, node distance=2.5cm]  {};
 \node[] (Q2) [above of =Q, node distance=1.01cm]  {};
\draw plot[smooth, tension=0.9] coordinates {(X1) (Y1) (P1) (X2)};
\path (X) edge (P);
\path (P) edge (Y);

 \node[] (Y) at (-1.5,-1.8)  {};   
 \node[] (P) at (-0.5,-0.6)  {};
 \node[] (Q) at (-0.166,-.2)  {};
 \node[] (Y1) [above of =Y, node distance=1.3cm]  {};   
 \node[] (P1) [above of =P, node distance=1.9cm]  {};
 \node[] (P2) [above of =P, node distance=1.02cm]  {};
 \node[] (Q1) [above of =Q, node distance=2.5cm]  {};
 \node[] (Q2) [above of =Q, node distance=1cm]  {};
\draw plot[smooth, tension=1.3] coordinates {(X1) (Y1) (P1) (X2)};
\path (X) edge (P);
\path (P) edge (Y);


\begin{scope}[xshift=16cm]

 \node[label={-60:\footnotesize $\mathcal{X}^B_s$}] (X) at (0,0)  {};
 \node[label={0:\footnotesize $\mathcal{X}^B_s$}] (X1) [above of =X, node distance=1cm]  {};
 \node[label={90:\footnotesize $\mathcal{X}^B_s$}] (X2) [above of =X,node distance= 3.2cm]  {};
 \node[] (Y) at (3,0)  {};   
 \node[] (P) at (1,0)  {};
 \node[] (Q) at (0.33,0)  {};
 \node[] (Y1) [above of =Y, node distance=2.9cm]  {};   
 \node[] (P1) [above of =P, node distance=2.3cm]  {};
 \node[] (P2) [above of =P, node distance=3.15cm]  {};
 \node[] (Q1) [above of =Q, node distance=1.6cm]  {};
 \node[] (Q2) [above of =Q, node distance=3.185cm]  {};
 \draw plot[smooth, tension=0.9] coordinates {(X1) (P1) (Y1) (X2)};
\path (X) edge (P);
\path (P) edge (Y);

 \node[] (Y) at (-2.1,0.6)  {};   
 \node[] (P) at (-0.7,0.2)  {};
 \node[] (Q) at (-0.233,0.066)  {};
 \node[] (Y1) [above of =Y, node distance=2.9cm]  {};   
 \node[] (P1) [above of =P, node distance=2.3cm]  {};
 \node[] (P2) [above of =P, node distance=3.15cm]  {};
 \node[] (Q1) [above of =Q, node distance=1.6cm]  {};
 \node[] (Q2) [above of =Q, node distance=3.19cm]  {};
\draw plot[smooth, tension=0.9] coordinates {(X1) (P1) (Y1) (X2)};
\path (X) edge (P);
\path (P) edge (Y);

 \node[] (Y) at (-1.5,-1.8)  {};   
 \node[] (P) at (-0.5,-0.6)  {};
 \node[] (Q) at (-0.166,-.2)  {};
 \node[] (Y1) [above of =Y, node distance=2.9cm]  {};   
 \node[] (P1) [above of =P, node distance=2.3cm]  {};
 \node[] (P2) [above of =P, node distance=3.15cm]  {};
 \node[] (Q1) [above of =Q, node distance=1.6cm]  {};
 \node[] (Q2) [above of =Q, node distance=3.18cm]  {};
\draw plot[smooth, tension=1] coordinates {(X1) (P1) (Y1) (X2)};
\path (X) edge (P);
\path (P) edge (Y);

\end{scope}
\end{tikzpicture}
\end{center}
\caption{The stable model $\mathfrak{T}_{\mathfrak{p}}$} \label{picTp}
\end{figure}

\par We now wish to increase the level structure at $\mathfrak{p}$. Restricting ourselves to the setting $B = \mathrm{M}_2(\Q)$, we now show how to obtain a semi-stable model for $U_p$ on $X_0(Np)$ from the work of Edixhoven \cite{Edi90}. To increase the number of cases in which $U_p$ is hyperbolic, and hence has a unique stable model, we puncture the modular curves at the cusps and treat $U_p$ as a punctured correspondence.
\begin{thmx}\label{UpB}
Let $p\geq 5$ coprime to $N$. Then the correspondence $U_p$ on $X_0(Np)$ has a skeletal semi-stable model $\mathfrak{U}_p$ over $W[\varpi]$ where $\varpi^{(p^2-1)/2} = p$, whose skeleton is depicted in Figure \ref{pic}. 
\begin{proof}
Denote the split Cartan curve appearing in $U_p$ by $X_{\sC}(N)$, then we have an isomorphism $X_{\sC}(N) \simeq X_0(p^2N)$. We use \cite[Proposition 10.4.6]{Liu02} to ensure that the normalisation $\mathcal{X}$ of the base change to $W[\varpi]$ of the blow-up the Katz--Mazur model $\mathcal{X}_{\sC}(N)$ at the supersingular points is semi-stable. Now extend scalars to $\mathcal{O}_{\C_p}$. We obtain exactly one additional component $Z$ on $\mathcal{X}_s$ for every supersingular point. We call $Z_i$ the component corresponding to the image under $\pi_i$ of the Type-II point corresponding to $Z$. As $Z$ is the unique pre-image of $Z_i$ under $\pi_i$, we obtain a skeletal semi-stable model $\mathfrak{U}_p$ after adding these components to the Katz--Mazur models of $X_0(Np)$, which is defined over $W[\varpi]$ by faithfully flat descent as in \ref{1Des}.

\par To determine the edge lengths, note that the expansion factors equal the ramification indices of $\pi_1$ and $\pi_2$ of the corresponding points on the ordinary components. As $\pi_1$ and $\pi_2$ restrict to the identity and Frobenius maps on the inner ordinary components, and are totally ramified above supersingular points on $X^{\pm}$, we obtain the expansion factors $1,p$ and $(p-1)/2$ as in Figure \ref{pic}, which uses the notation of Edixhoven for the ordinary components using the isomorphism $X_{\sC}(N) \simeq X_0(Np^2)$. The numbers marking the edges are the expansion factors, which we omitted whenever they are $1$. 
\end{proof}
\end{thmx}

\begin{figure}[htbp!]
\begin{center}
\begin{tikzpicture}[-,scale=.5]
 \path (2,3.5) edge[->,>=stealth'] node[above right,draw=none] {\footnotesize $\pi_2$} (3,2);
 \path (-7,3.5) edge[->,>=stealth'] node[above left,draw=none] {\footnotesize $\pi_1$} (-8,2);
 \node[label={=}] at (-2.5,4.6) {};
 \node[label={90:\tiny $\infty$}] at (.8,-.15) {};
 \node[label={90:\tiny $0$}] at (7.2,-.15) {};
\begin{scope}[xscale=-1,xshift=5cm]
 \node[label={90:\tiny $0$}] at (.8,-.15) {};
 \node[label={90:\tiny $\infty$}] at (7.2,-.15) {};
\end{scope}

\tikzstyle{every node}=[draw,circle,fill=black,minimum size=2pt,inner sep=0pt]
 \node[] (XF) at (2,0)  {};
 \node[]   (XV) at (6,0)  {};   
 \node[label={90:\scriptsize $Z_2$}]    (1) at (3,1)  {};
 \node[label={90:\scriptsize $Z_2$}]    (2) at (3,0)  {};
 \node[label={90:\scriptsize $Z_2$}]    (3) at (3,-1) {};
\draw plot[smooth cycle, tension=0.7] coordinates {(XF) (1) (XV) (3)};
\path (XF) edge (XV);
\path (XF) edge (1.5,0);
\path (1.5,0) edge[dashed] (0.8,0);
\path (XV) edge (6.5,0);
\path (6.5,0) edge[dashed] (7.2,0);

 \node[]                               (XFF) at (0,5)  {};
 \node[]                               (XVV) at (4,5)  {}; 
 \node[label={0:\scriptsize $X^+$}]  (X+) at (0,8)  {};
 \node[label={0:\scriptsize $X^-$}]  (X-) at (0,2)  {};  
 \node[label={70:\scriptsize $Z$}]     (11) at (1,6)  {};
 \node[label={70:\scriptsize $Z$}]     (22) at (1,5)  {};
 \node[label={70:\scriptsize $Z$}]     (33) at (1,4)  {};
\draw plot[smooth cycle, tension=0.7] coordinates {(XFF) (11) (XVV) (33)};
\path (XVV) edge (XFF);

 \path (X+) edge[bend right=20] (11);
 \path (X+) edge[bend right=15] (22);
 \path (X+) edge[bend right=15] (33);  
 \path (X-) edge[-,line width=3pt,draw=white,bend left=15] (22);
 \path (X-) edge[-,line width=3pt,draw=white,bend left=15] (11);
 \path (X-) edge[bend left=20] (33);
 \path (X-) edge[bend left=15] (22);
 \path (X-) edge[bend left=15] (11); 
 \path (XFF) edge (-.5,5) ;
 \path (-.5,5) edge[dashed] (-1.2,5);
 
 \path (XVV) edge (4.5,5)
 			 edge[bend left] (4.5,5.3)
 			 edge[bend right] (4.5,4.7);
 \path (4.5,5) edge[dashed] (5.2,5);
 \path (4.5,5.3) edge[dashed] (5.2,5.3);
 \path (4.5,4.7) edge[dashed] (5.2,4.7);
\begin{scope}[xshift=-9cm]
 \path (0,5) edge (-.5,5);
 \path (-.5,5) edge[dashed] (-1.2,5);
 \path (4,5) edge (4.5,5)
 			 edge[bend left]  (4.5,5.3) 
 			 edge[bend right] (4.5,4.7);
 \path (4.5,5) edge[dashed]   (5.2,5);
 \path (4.5,5.3) edge[dashed] (5.2,5.3);
 \path (4.5,4.7) edge[dashed] (5.2,4.7);
\end{scope}

 \node[white,label={0:\tiny $\mathrm{Fr}_p$}] at (1.4,5.5) {};
\path (1.2,5.95) edge[<->,>=stealth',bend left=60] (1.2,5.05);

\begin{scope}[xscale=-1,xshift=5cm]
 \node[] (XF) at (2,0)  {};
 \node[]   (XV) at (6,0)  {};   
 \node[label={90:\scriptsize $Z_1$}]    (1) at (3,1)  {};
 \node[label={90:\scriptsize $Z_1$}]    (2) at (3,0)  {};
 \node[label={90:\scriptsize $Z_1$}]    (3) at (3,-1) {};
\draw plot[smooth cycle, tension=0.7] coordinates {(XF) (1) (XV) (3)};
\path (XF) edge (XV);
\path (XF) edge (1.5,0);
\path (1.5,0) edge[dashed] (0.8,0);
\path (XV) edge (6.5,0);
\path (6.5,0) edge[dashed] (7.2,0);

 \node[]                             (XFF) at (0,5)  {};
 \node[]                             (XVV) at (4,5)  {}; 
 \node[label={180:\scriptsize $X^+$}](X+)  at (0,8)  {};
 \node[label={180:\scriptsize $X^-$}](X-)  at (0,2)  {};  
 \node[label={110:\scriptsize $Z$}]  (11)  at (1,6)  {};
 \node[label={110:\scriptsize $Z$}]  (22)  at (1,5)  {};
 \node[label={110:\scriptsize $Z$}]  (33)  at (1,4)  {};

\draw plot[smooth cycle, tension=0.7] coordinates {(XFF) (11) (XVV) (33)};
\path (XVV) edge (XFF);

 \path (X+) edge[bend right=20] (11);
 \path (X+) edge[bend right=15] (22);
 \path (X+) edge[bend right=15] (33);  
 \path (X-) edge[-,line width=3pt,draw=white,bend left=15] (22);
 \path (X-) edge[-,line width=3pt,draw=white,bend left=15] (11);
 \path (X-) edge[bend left=20] (33);
 \path (X-) edge[bend left=15] (22);
 \path (X-) edge[bend left=15] (11);  
\end{scope} 

\begin{scope}
 \tikzstyle{every node}=[circle,fill=white,minimum size=8pt,inner sep=0pt]
 \node at (2.5,4.3) {\tiny $p$};
 \node at (2.5,5)   {\tiny $p$};
 \node at (2.5,5.7) {\tiny $p$};
 \node at (-.5,7)   {\tiny $\frac{p-1}{2}$};
 \node at (-.5,3)   {\tiny $\frac{p-1}{2}$};
\begin{scope}[xscale=-1,xshift=5cm]
 \node at (2.5,4.3)   {\tiny $p$};
 \node at (2.5,5)   {\tiny $p$};
 \node at (2.5,5.7) {\tiny $p$};
 \node at (-.5,7)   {\tiny $\frac{p-1}{2}$};
 \node at (-.5,3)   {\tiny $\frac{p-1}{2}$};
\end{scope}
\end{scope}
\end{tikzpicture}
\end{center}
\caption{The skeleton of $\mathfrak{U}_p$} \label{pic}
\end{figure}

\par As a consequence of Theorem \ref{UpB}, we obtain that $U_p$ acts on the toric part of $\HH^1_{\et}(X_0(Np)_{\overline{\Q}_p},\Q_l)$ as Frobenius. This is \cite[Proposition 3.8]{Rib90}. We also recover the action on the N\'eron component group $\Phi$. 

\subsection{Hecke operators at the infinite level. }\label{4Inf}As we increment the level at $\mathfrak{p}$ of our Hecke operators, it becomes increasingly difficult to find semi-stable models. In view of the recent work of Weinstein \cite{Wei12}, it seems prudent to work at the infinite level instead. We opted for the language of adic spaces, partly because of the more streamlined treatment of the main definitions due to the presence of Type-V points, but also because we expect the proof of our main theorem \ref{ss} to remain valid in the setting of perfectoid spaces. If one can give a satisfactory description of the harmonic maps constituting $U_p$ at the infinite level, this would likely reveal interesting information. We hope to return to this problem in the near future. 

\subsection{Too supersingular elliptic curves. }\label{Moduli}Coleman \cite{Col05} gives a moduli interpretation of the components in the special fibre of Edixhoven's semi-stable model $\mathcal{X}_0(Np^2)$ in tame level $N=1$. As an immediate consequence of the determination of the stable model $\mathfrak{U}_p$ given above, we now generalise this theorem to general $N$. Coleman's argument is roughly equivalent to ours, and should extend to any tame level $N$.

\par \noindent \textbf{Terminology. }An elliptic curve over $\C_p$ is said to be \textit{too supersingular} if it does not have a canonical subgroup as in \cite{Kat73} and \cite{Lub79}. It is \textit{nearly too supersingular} if it is $p$-isogenous to such a curve. 

\par The ordinary components of $\mathfrak{U}_p$ arise as strict transforms of components in the Katz--Mazur models, or normalisations thereof, and hence carry a natural moduli interpretation. The following corollary of Theorem \ref{UpB} provides a moduli interpretation for the supersingular components. 

\begin{thmx}\label{Mod}Resume the notation of Figure \ref{pic}. The points in $\mathrm{sp}^{-1}(Z_1^{\sm})$ and $\mathrm{sp}^{-1}(Z_2^{\sm})$ parametrise too supersingular curves, resp. nearly too supersingular curves, together with $\Gamma_0(N)$-level structure.
\begin{proof}
For an open annulus $A$ with outer radius $1$ and inner radius $r$, the open subannulus $A_s := \{|.|_x : 1 \geq |T|_x \geq s\}$ for $s \geq r$ is independent of the chosen parameter $T$ for $A$. Letting $E_{p-1}$ be a global parameter for the supersingular annuli in $X_0(Np)$, lifting the Hasse invariant, we see that the smooth points of $Z_1$ and $Z_2$ classify those elliptic curves over $\mathcal{O}_{\C_p}$ with Hasse invariants $p/(p+1)$ and $1/(p+1)$ respectively. The claim follows from the work of Katz \cite{Kat73} and Lubin \cite{Lub79}. 
\end{proof}
\end{thmx}
\section{Overconvergent \texorpdfstring{$p$}{p}-adic modular forms}\label{5Chap}We now discuss some consequences for the theory of overconvergent modular forms. We generalise the work of Goren--Kassaei \cite{GK06} on canonical subgroups for a general class of curves, including Shimura curves over totally real fields, and discuss the work of Kassaei \cite{Kas09} on analytic continuation.

\par The $p$-adic variation of modular forms is described by the theory of \textit{overconvergent $p$-adic modular forms}, as developed primarily by Katz \cite{Kat73}, Coleman \cite{Col97} and Coleman--Mazur \cite{CM98}. We will now discuss two of the main results of the theory: canonical subgroups and analytic continuation. Over $\mathrm{GL}(2)_{\Q}$ these results were proved using the moduli interpretation of modular curves. The work of Goren--Kassaei \cite{GK06} and Kassaei \cite{Kas09} avoids moduli interpretations, and focusses instead on the geometry of the curves involved. Their approach has proved fruitful for subsequent generalisations to other reductive groups, most notably the theory of overconvergent Hilbert modular forms \cite{GK09,GK12}.

\subsection{Canonical subgroups for general curves. }\label{Canon}Let $K$ be a field complete with respect to a non-trivial non-Archimedean valuation $|\cdot |$ of rank $1$, with valuation ring $R$. We generalise \cite[Theorem 3.9]{GK06}.
\begin{thmx}\label{Can}Let $g:\mathcal{X} \rightarrow \mathcal{Y}$ be a finite flat morphism between semi-stable curves over $R$, and let $\mathfrak{g}:\mathcal{X}_1 \rightarrow \mathcal{Y}_1$ be the stable hull of $\C$ with respect to $(\mathcal{X},\mathcal{Y})$, see Theorem \ref{ssskel}. Assume that $g$ is an isomorphism when restricted to a component $Z$ of $\mathcal{X}_{s}$. Then $g$ is an isomorphism on $\mathrm{sp}_{\mathcal{X}_1}^{-1}(Z)$.
\begin{proof}
Every edge adjacent to $Z$ in the skeleton of $\mathfrak{g}$ must have expansion factor $1$, as this is the ramification index of the corresponding point on $Z$. We obtain the following picture:
\begin{center}
\begin{figure}[htbp!]
\begin{tikzpicture}[-,scale=.3,node distance=1.8cm]
 \path (-2,0) edge[->,>=stealth'] node[above,draw=none] {$\mathfrak{g}$} (1,0);

\tikzstyle{every node}=[draw,circle,fill=black,minimum size=1.7pt,inner sep=0pt]

 \node[label={180:\mbox{\scriptsize $\simeq Z$}}] (X) at (5,0)  {};
 \node[] (E1) at (5,3)  {};   
 \node[] (E2) at (4,-3) {}; 
 \node[] (E3) at (7,-5) {};  
 \node[] (IN) at (6,2) {};  
 \node[label={180:\mbox{\scriptsize $Z$}}] (X1) at (-10,0)  {};
 \node[] (Y1) at (-9.4,3)  {};
 \node[] (Y2) at (-11.2,-3){};
 \node[] (Y3) at (-8.4,-5) {};
 \node[] (Z1) at (-7,0) {};
 \node[] (Z2) at (-14,0){};
 \node[] (Z3) at (-5,0) {};
 \node[] (N)  at (-6,2) {};
 \node[label={180:$\mathrm{sp}_{\mathcal{X}_1}^{-1}(Z)$},fill=none,minimum size=0pt] (XZ) at (-12,2.5) {};

\draw plot[smooth, tension=0.9] coordinates {(X1) (-9.6,3) (-8,0.6) (Z1)};
\draw plot[smooth, tension=0.9] coordinates {(X1) (Y2) (-12.8,-1) (Z2)};
\draw plot[smooth, tension=0.9] coordinates {(X1) (Y3) (-7,-1.6) (Z3)};
\draw plot[smooth, tension=0.9] coordinates {(Z1) (N) (Z3)};
\draw plot[smooth, tension=0.9] coordinates {(X) (E1) (IN) (X)};
\draw [dashed] plot[smooth, tension=0.4] coordinates {(Y3) (Y1) (-11,2) (Y2) (-9.4,-5) (Y3)};
\path (X)  edge (E2) edge (E3);
\path (N)  edge[bend right=15] (Y1);
\path (-10.6,2)  edge[bend right=30] (-12,2.5); 

\begin{scope}
 \tikzstyle{every node}=[draw,circle,fill=white,minimum size=7pt,inner sep=0pt]
 \node at (-10,1.5) {\tiny $1$};
 \node at (-10.4,-1.8) {\tiny $1$};
 \node at (-9.3,-3) {\tiny $1$};
\end{scope}
\end{tikzpicture}
\end{figure}
\end{center}
Following Goren--Kassaei \cite[Proposition 3.1]{GK06}, the section $\mathfrak{g}(Z) \rightarrow Z$ gives rise to a section $s:\mathrm{sp}^{-1} (\mathfrak{g}(Z)^{\sm}) \rightarrow X^{\ad}$ on the tube of the smooth locus $\mathfrak{g}(Z)^{\sm}$, as $R$ is complete and therefore Henselian. By Coleman--Gouv\^ea--Jochnowitz \cite[Lemma 6]{CGJ95}, we may extend $s$ to a section $s^{\dagger}:U \rightarrow X^{\ad}$ for $U \subset Y^{\ad}$ an open subspace strictly containing the domain of $s$. As the expansion factors of the edges in the skeleton of $\mathcal{X}_1$ adjacent to $Z$ are all $1$, we see that $g$ defines an isomorphism between the wide open $V = \mathrm{sp}^{-1}_{\mathcal{X}_1}(Z^{\sing})$ and its image. Its inverse agrees with $s^{\dagger}$ on $U \cap \mathfrak{g}(V)$, and hence we may glue it together with $s^{\dagger}$ to produce the desired section.
\end{proof}
\end{thmx}

\subsection{Remark. }Taking the union of the open regions $\mathrm{sp}_{\mathcal{X}_1}^{-1}(Z)$ over all components $Z$ on which $g$ induces an isomorphism, we easily see that the resulting section is maximal in the obvious sense. The value of Theorem \ref{Can} is that it identifies this maximal section explicitly, in terms of the stable hull of $g$. A special case was treated in Goren--Kassaei \cite[Theorem 3.9]{GK06}, where $\mathcal{X}_s$ was assumed to be composed of two components, intersecting transversally, $\mathcal{Y}$ was assumed to have good reduction, and $K$ was a finite extension of $\Q_p$. These assumptions are satisfied in the case of Shimura curves $X^B_0(\mathfrak{p}) \rightarrow X^B$. It is straightforward to check that $g^{\ad}$ maps $\mathrm{sp}_{\mathcal{X}_1}^{-1}(Z)$ isomorphically to the adic space attached to the rigid subspace $\mathfrak{Y}_{\rig}[0,e/(e+1))$ defined by Goren--Kassaei \cite[Theorem A]{GK06}, which recovers their result.

\subsection{Analytic continuation. }One crucial ingredient in the theory of overconvergent modular forms is Coleman's classicality criterion. It guarantees that an overconvergent modular eigenform of weight $k$ is in fact a classical modular form if the $p$-adic valuation of its $U_p$-eigenvalue is less than $k-1$. Coleman \cite{Col96} initially proved this criterion cohomologically, an approach which is taken by Johansson \cite{Joh13} but which seems difficult to generalise further. Buzzard \cite{Buz03} and Kassaei \cite{Kas06} develop an approach which involves $p$-adically continuing an overconvergent modular form to the entire supersingular locus, and gluing it together with an appropriately constructed form on a complementary region arising from the Atkin operator. Kassaei \cite{Kas09} proves analytic continuation results following this strategy, but only referring to the underlying geometry of the $U_p$ operator. His results may naturally be viewed in the framework of semi-stable models for correspondences developed in this paper.

\subsection{Higher dimensions. }\label{Dim}It is natural to wonder whether Theorem \ref{ss} generalises to correspondences between varieties of arbitrary dimension. If $X$ is a smooth and proper variety over $K$, one can attach to any semi-stable formal model $\mathfrak{X}$ a skeleton of $X^{\ad}$. Goren--Kassaei \cite{GK09,GK12} prove the existence of canonical subgroups for Hilbert modular varieties using methods close in spirit to our skeletal approach, which suggests that some generalisation for correspondences in higher dimensions might exist. 

\bibliographystyle{alpha}
\bibliography{/Users/vonk/Documents/Wiskunde/TeX/References}
\end{document}